\documentclass[12pt]{amsart}

\usepackage{fullpage}

\usepackage{amsmath}
\usepackage{amsfonts}
\usepackage{amssymb}
\usepackage{amsthm}
\usepackage{mathtools} 
\usepackage{latexsym}

\usepackage{enumerate}
\usepackage{cancel}
\usepackage{cases}
\usepackage{empheq}
\usepackage{multicol}

\usepackage{wrapfig}

\theoremstyle{plain}

\theoremstyle{definition}

\theoremstyle{remark}

\numberwithin{equation}{section} 
\numberwithin{figure}{section}   

\usepackage{cleveref}

\newcommand{\field}[1]{\mathbb{#1}}
\newcommand{\nN}{\field{N}}

\newcommand{\nT}{\mathbb T}

\newcommand{\set}[1]{\left\{#1\right\}}

\newcommand{\norm}[1]{\left\|#1\right\|}

\newcommand{\grad}{\nabla}
\renewcommand{\div}{\nabla \cdot} 
\newcommand{\lap}{\triangle} 

\usepackage{placeins} 


\author{Trenton Franz}
\address[Trenton Franz]{School of Natural Resources, 
                University of Nebraska--Lincoln,
        Lincoln, NE 68588-0130, USA}
\email[Trenton Franz]{trenton.franz@unl.edu}
\author{Adam Larios}
\address[Adam Larios]{Department of Mathematics, 
                University of Nebraska--Lincoln,
        Lincoln, NE 68588-0130, USA}
\email[Adam Larios]{alarios@unl.edu}
\author{Collin Victor}
\address[Collin Victor]{Department of Mathematics, 
                University of Nebraska--Lincoln,
        Lincoln, NE 68588-0130, USA}
\email[Collin Victor]{collin.victor@huskers.unl.edu}

\title{The Bleeps, the Sweeps, and the Creeps:\\Convergence Rates for Dynamic Observer Patterns via Data Assimilation for the 2D Navier-Stokes Equations}
\date{\today}

\begin{document}

\begin{abstract}
We adapt a continuous data assimilation scheme, known as the Azouani-Olson-Titi (AOT) algorithm, to the case of moving observers for the 2D incompressible Navier-Stokes equations.  We propose and test computationally several movement patterns (which we refer to as ``the bleeps, the sweeps and the creeps''), as well as Lagrangian motion and combinations of these patterns, in comparison with static (i.e. non-moving) observers.  In several cases, order-of-magnitude improvements in terms of the time-to-convergence are observed.  We end with a discussion of possible applications to real-world data collection strategies that may lead to substantial improvements in predictive capabilities.
\end{abstract}

\maketitle

\thispagestyle{empty}

\noindent
{\small
``\textit{I'm having trouble with the radar, sir. [...] I've lost the bleeps, I've lost the sweeps, and I've lost the creeps.}'' --Michael Winslow as Radar Technician in \textit{Spaceballs} \cite{spaceballs}.
}

\section{Introduction}
\noindent
One difficulty in the simulation of many physically interesting dynamical systems is that initial data is often incomplete or inaccurate. However, one can circumvent this difficulty by using the tools of \textit{data assimilation}, which is a class of techniques used to increase accuracy by combining time-dependent observational data together with a physical model, driving simulations to converge to the ``true'' solution.  In this work, we concentrate on a simple yet powerful data assimilation tool known as the Azouani-Olson-Titi (AOT) algorithm, since it readily lends itself to the case of \textit{dynamic observation points}, which is the major focus of this manuscript.  In a previous work \cite{Larios_Victor_2019}, two of the authors first proposed the idea of dynamic observers for the AOT algorithm, and showed that in a simplified 1D model (the Allen-Cahn reaction-diffusion equations), moving the observation points in time resulted in order-of-magnitude increases in the time to convergence to machine precision error.  In the present work, we extend the idea of continuous data assimilation via moving observers to the 2D Navier-Stokes equations, and show similar gains in convergence rates.  Moreover, we propose and investigate several movement schemes for observers, and compare these to the case of static (i.e., non-moving) observers.   During the completion of this manuscript, we learned of another work \cite{Biswas_Bradshaw_Jolly_2020}, which shows convergence of the AOT algorithm for the 2D Navier-Stokes equations with observers in a moving square patch.  However, the purpose of the present work is to examine the effectiveness of several different, arguably more realistic, movement patterns with the aim of real-world data collection strategies.

Since its inception in \cite{Azouani_Olson_Titi_2014,Azouani_Titi_2014}, the AOT algorithm has been the subject of much recent study in both analytical studies  
\cite{
Albanez_Nussenzveig_Lopes_Titi_2016,
Bessaih_Olson_Titi_2015,
Biswas_Bradshaw_Jolly_2020,
Biswas_Foias_Mondaini_Titi_2018downscaling,
Biswas_Martinez_2017,
Biswas_Price_2020_AOT3D,
Carlson_Hudson_Larios_2018,
Carlson_Larios_2021_sens,
Celik_Olson_Titi_2019,
Chen_Li_Lunasin_2021,
Diegel_Rebholz_2021,
Du_Shiue2021,
Farhat_Jolly_Titi_2015,
Farhat_Lunasin_Titi_2016abridged,
Farhat_Lunasin_Titi_2016benard,
Farhat_Lunasin_Titi_2016_Charney,
Farhat_Lunasin_Titi_2017_Horizontal,
Foias_Mondaini_Titi_2016,
Foyash_Dzholli_Kravchenko_Titi_2014,
GarciaArchilla_Novo_Titi_2018,
GarciaArchilla_Novo_2020,
Gardner_Larios_Rebholz_Vargun_Zerfas_2020_VVDA,
Ibdah_Mondaini_Titi_2018uniform,
Jolly_Martinez_Olson_Titi_2018_blurred_SQG,
Jolly_Martinez_Titi_2017,
Jolly_Sadigov_Titi_2015,
Larios_Pei_2018_NSV_DA,
Larios_Victor_2021_chiVsdelta2D,
Markowich_Titi_Trabelsi_2016_Darcy,
Mondaini_Titi_2018_SIAM_NA,
Pachev_Whitehead_McQuarrie_2021concurrent,
Pei_2019,
Rebholz_Zerfas_2018_alg_nudge,
Zerfas_Rebholz_Schneier_Iliescu_2019} 
and computational studies 
\cite{
Altaf_Titi_Knio_Zhao_Mc_Cabe_Hoteit_2015,
Carlson_Hudson_Larios_Martinez_Ng_Whitehead_2021,
Carlson_VanRoekel_Petersen_Godinez_Larios_2021,
DiLeoni_Clark_Mazzino_Biferale_2018_inferring,
Desamsetti_Dasari_Langodan_Knio_Hoteit_Titi_2019_WRF,
Gesho_Olson_Titi_2015,
Larios_Pei_2017_KSE_DA_NL,
Larios_Rebholz_Zerfas_2018,
Larios_Victor_2019,
Lunasin_Titi_2015}.

In this work we investigate the AOT algorithm as applied to the 2D Navier-Stokes equations (NSE). Writing the NSE in an abstract form\footnote{After solving for the pressure, one can write $F(u)=\nu\lap u-u\cdot\nabla u - \nabla\lap^{-1}\nabla\cdot(u\otimes u)$ (see, e.g., \cite{Constantin_Foias_1988,Temam_1995_Fun_Anal} for details), but this specific form is not especially relevant for present discussion.} we have the following equation  
\begin{equation}\label{eq_dyn_sys}
    \frac{d}{dt}u = F(u).\\
\end{equation}
Note that in \eqref{eq_dyn_sys}, we do not know the initial condition for this system.  However, we assume that we can observe the solution at certain spatial locations, which may be sparsely distributed, say with some characteristic spatial distance between observations $h>0$. Using these measurements, the AOT algorithm \cite{Azouani_Olson_Titi_2014,Azouani_Titi_2014} constructs the system
\begin{align}\label{eq_dyn_sys_AOT}\begin{cases}
    \frac{d}{dt}v = F(v) + \mu (I_h(u)-I_h(v))\\
    v(0) = v_0
    \end{cases}
\end{align}
where $I_h(u)$ is a suitable interpolation of the observed values of a strong solution $u$ of \eqref{eq_dyn_sys_AOT}, and $\mu>0$ is a constant, sometimes called the ``nudging parameter.'' It was proven in \cite{Azouani_Olson_Titi_2014} that for sufficiently large $\mu$ and sufficiently small $h$ that  $\norm{u(t)-v(t)}_{L^2}\to 0$ and $\norm{\nabla u(t)-\nabla v(t)}_{L^2}\to 0$ exponentially fast in time as $t\to\infty$ for \textit{any} choice of $v_0 \in L^2$ such that $v_0$ is divergence-free. Note that in many works including the present one, for simplicity data is assumed to be observed continuously in time and without noise.  However, similar exponential convergence results for discrete-in-time observations were proven in \cite{Foias_Mondaini_Titi_2016}, and for noisy data and stochastic forcing in \cite{Bessaih_Olson_Titi_2015,Carlson_Hudson_Larios_Martinez_Ng_Whitehead_2021}.

An important aspect to consider when capturing physical phenomena via data assimilation is where to place observers. Observers can be thought of as physical objects such as weather stations in the case of predicting the weather. These observers may require financial investment and manpower to place, thus it is crucial to minimize the number of observers required for convergence of the method. Computational tests in the case of 2D NSE indicated exponential convergence in time using a static uniform grid of observers \cite{Gesho_Olson_Titi_2015}. The rate of convergence can typically be sped up by adding more observers, however computational studies \cite{Biswas_Bradshaw_Jolly2020,Larios_Victor_2019} seem to suggest that the rate of convergence can be improved by moving existing observers instead.
In a previous work we found in the case of the 1D Allen-Cahn equations that the number of observers can be drastically reduced by replacing the static uniform grid with a mobile cluster of fine length scale \cite{Larios_Victor_2019}.
Recently in \cite{Biswas_Bradshaw_Jolly2020} it was shown computationally that mobile local data assimilation in the case of 2D NSE achieved fast convergence rates when compared to a non-mobile local data assimilation algorithm. There mobile local data assimilation was investigated in the single regime of continuously moving a square subdomain around the domain. In this work we investigate computationally the effect of movement strategies for observers on convergence rates for 2D NSE for both local and global data assimilation.

The effect of sweeps has been recently observed in the reduction of skill in weather forecasting models due to the elimination of cross Atlantic air travel. For example, it was shown in \cite{Chen_2020covid} the 50-75\% reduction in air flights in 2020 has lead to a decrease of $2^\circ$C in air temperature forecasts over Greenland and Siberia and a sustained reduction in flights will lead to further degradation of the air pressure and wind speed forecast skill.

 
In the present work, we focus on the 2D incompressible Navier-Stokes equations, given by:
\begin{subequations}\label{NSE}
\begin{align}
    u_t + u\cdot \grad u - \nu \lap u &= \grad p + f,\\
    \div u &= 0.
\end{align}
\end{subequations}
Here $u(x,t) = (u_1(x,t),u_2(x,t))$ is the velocity, $\nu>0$ is the viscosity, $p$ is the pressure and $f$ is a forcing term. Applying the continuous data assimilation algorithm \eqref{eq_dyn_sys_AOT} to this system yields the following equations:
\begin{subequations}\label{NSE-v}
\begin{align}
    v_t + v\cdot \grad v - \nu \lap v &= \grad p + f + \mu (I_{h,t}(u)-I_{h,t}(v)),\\
    \div v &= 0.
\end{align}
\end{subequations}
Here, $I_{h,t}$ is a linear operator representing nodal piecewise-linear interpolation, and $\mu>0$ is a algorithmic parameter often referred to as the ``nudging'' parameter.  In the static case, more general interpolation operators were considered in  \cite{Azouani_Olson_Titi_2014}, but for the sake of simplicity, we only consider the case of piecewise-linear interpolation.   
Note that, in contrast to the usual AOT algorithm for static observers, here, we allow the observer locations to depend on time, as denoted by the subscript $t$ in $I_{h,t}$.

This paper is laid out as follows.  
In Section 2, we describe the various observer movement patterns that we consider.  
In Section 3, we describe our computational methods and results, but to keep the discussion concise, we only refer to the graphics, which are in the appendix. 
In Section 4, we conclude with a comparative discussion on the merits of the various observer movement patterns, and contrast them with the results of the static observer vase.

\section{Movement Patterns}
For this study we computationally examined different methods of moving observers dynamically.  We group these into three general movement patterns that we call ``The Bleeps'', ``The Sweeps'' and ``The Creeps''. In addition to these movement patterns we also investigate a fourth pattern given by observers as Lagrangian particles suspended in the fluid.


\subsection{The Bleeps}
``The Bleeps'' refers to the strategy of observing data at random locations that change over time. Initially observers are placed at random locations and at each time-step.  The old observers are then removed, and new observers are placed at new random locations. These random locations are determined by picking uniformly random integers $x_i \in \set{0,1,...,N}$ with locations given by $x = (x_1\cdot \Delta x - \pi ,x_2 \cdot \Delta y - \pi)$. The locations are chosen in this manner so that their coordinates align with the underlying spatial grid.

While we do not claim that there is a  physically-realistic scenario corresponding to this movement strategy, studying ``The Bleeps'' can give some idea of how non-static observations can influence convergence rates in the absence of any coherent spatial structure.  Thus, we examine this case as guideline or comparative test to shed light on the importance of domain coverage vs. certain geometric patterns in observer motion.

\subsection{The Sweeps}
``The Sweeps'' refers to the general strategy of moving a group observers at a constant velocity periodically over the domain. For this method we investigated three separate movement strategies referred to as thin, thick, and random sweeps.

``Thin-sweeps'' refers to the strategy of utilizing a large cluster of observers placed in a thin rectangular strip of the domain that moves continuously to the right at a constant velocity. The rectangular strip utilized by this strategy has length $x = a\Delta x$ and width $y = L_y = 2\pi$ and moves at constant velocity $b\frac{\Delta x}{\Delta t}$ for some $a,b\in \nN$. Note that the form of the velocity and $x$-dimension were chosen specifically so that the rectangular region would align with the underlying spatial grid at every timestep. In all of the simulations conducted we choose values $a\in \set{1,...,20}$ and additionally with $b=a$.
In simulations this rectangular box contains observers placed at each gridpoint along the spatial grid, meaning that the full solution of $u$ is observed locally in this rectangular region. Outside of this rectangle, $I_{h,t}(u-v) = 0$ for any given time $t$. Thus $I_{h,t}$ is given as follows:
\begin{align*}
    I_{h,t}(u-v) = (u-v)\chi_{R(t)}
\end{align*}
Here $R(t)$ refers to the sliding rectangular box of observers at the $t$. 

``Thick-sweeps'' refers to the strategy of utilizing a uniform grid of observers placed on a rectangular region that does not vary in size with the number of observers contained within. We fix the rectangular region to be a box containing the fourth of the domain given by $R = \set{ (x,y)\in \nT^2: -\pi \leq x \leq -\frac{\pi }{2} }$. This rectangular region moves to the right at the constant velocity $b\frac{\Delta x}{\Delta t}$. In all of the simulations whose results are presented here we used $b = 1$.
For certain velocities one can view this method as a uniform static grid with only a fourth of the observers active at any given time.
$I_{h,t}$ is given as follows:
\begin{align*}
    I_{h,t}(u-v) =  (I_{h}(u-v))\chi_{R(t)} 
\end{align*}

\begin{figure}[htp!]
	\includegraphics[width=\textwidth]{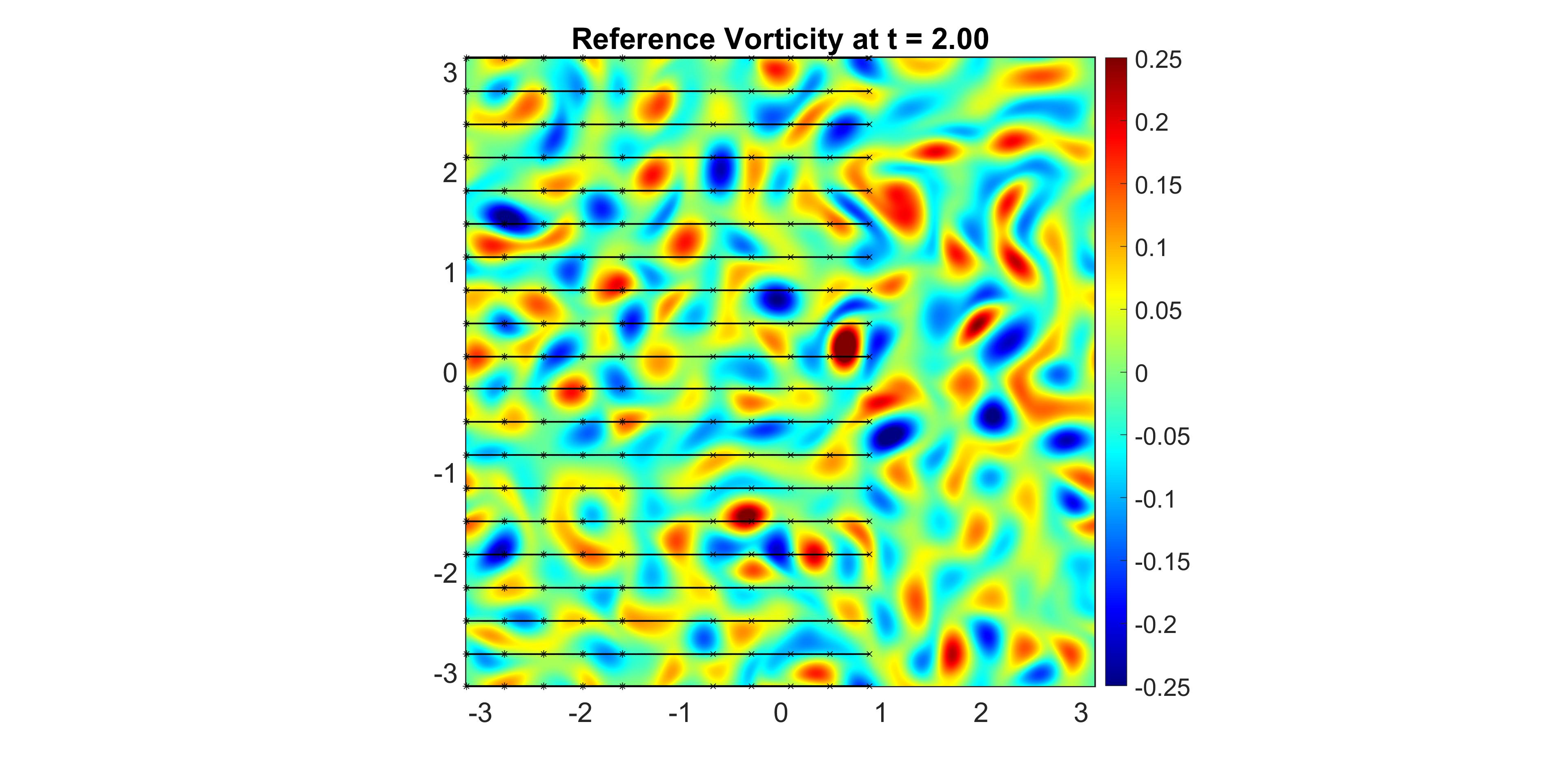}
	\caption{Graph displaying movement of observers using the thick sweeps movement strategy displayed over the vorticity of the reference solution at time $t=2$. Observers were initially placed at locations indicated by * and travel to the right along the black paths with ending locations marked by $\times$.}
	\label{fig:MovementBox}
\end{figure}

Random-sweeps refers to the strategy of placing observers at random locations with an additional random velocity.
A major difference between this and all of the previous methods is that the location of the observers are not restricted to the underlying spatial grid. This means that the value of $u-v$ must be interpolated from the spatial grid at each timestep in order obtain approximate values for the observer to observe. The velocities of each observer are given by vectors with each component having a random uniformly distributed value from the interval $(-1,1)$.
$I_{h,t}$ is given as follows:
\begin{align*}
    I_{h,t}(u-v) = I_{h,t} \circ J_{\Delta x,t}(u-v)
\end{align*}
Here $J$ is a linear interpolation operator applied to $u-v$ in order to get the approximate values for the observers with $I$ interpolating these observed values onto $\nT^2$.

\begin{figure}[htp!]
	\includegraphics[width=\textwidth]{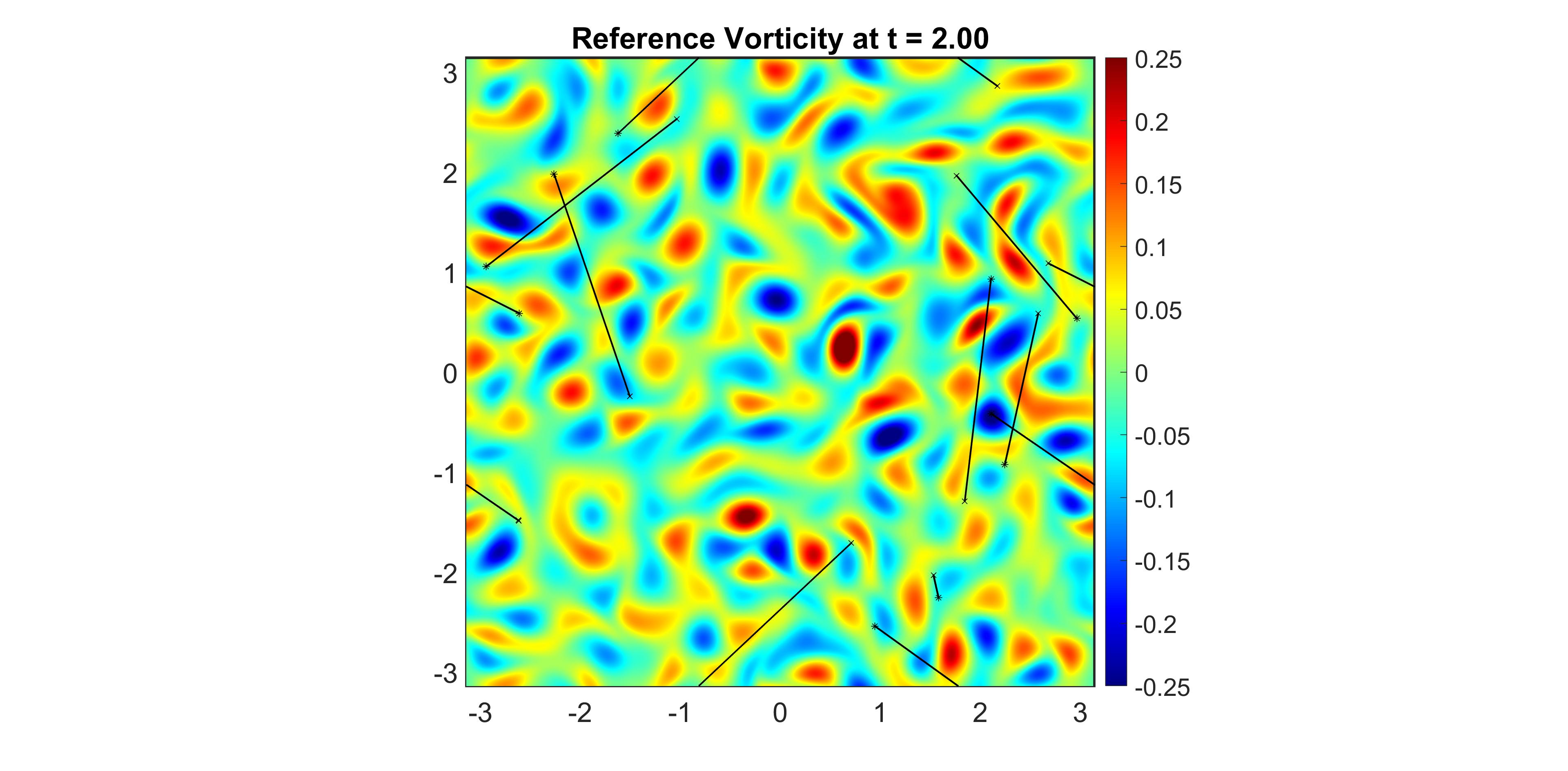}
	\caption{Graph displaying movement of observers using the random sweeps movement strategy displayed over the vorticity of the reference solution at time $t=2$. Observers were initially placed at 10 locations indicated by * and travel along random trajectories indicated by black paths with ending locations marked by $\times$.}
	\label{fig:MovementRandom}
\end{figure}

We study each of these sweep-movement strategies as each can be potentially realized in physical scenarios. The thin-sweep observers can be used to model observations given by, e.g., a laser sweeping through a domain in particle image velocimetry (PIV) applications, where a velocity profile is taken as an observation. Thick sweeps can be used to model observers on mobile platforms that move at a fixed velocity. Random sweeps can be thought of as modeling observers with a fixed velocity, e.g., airplanes or satellites.


\subsection{The Creeps}
``The creeps'' refers to the strategy of having observers move according to random walks (``creeping'' along). After a given amount of time, each observer will (uniformly) randomly pick an integer $j\in \set{1,2,3,4,5}$.
These integers determine the direction of motion $d(j)$ as follows:
\begin{align*}
    d(j) =\left(\begin{matrix*}x_1\Delta x\\x_2\Delta y \end{matrix*}\right) \
    \text{ where }
    \begin{cases}x_1=1, x_2=0,& \text{ if } j = 1,\\
    x_1=-1, x_2=0, & \text{ if } j=2,\\
    x_1=0, x_2=1, & \text{ if } j=3,\\
    x_1=0, x_2=-1, & \text{ if } j=4,\\
    x_1 = 0, x_2 = 0 & \text{ if } j=5.
    \end{cases}
\end{align*}
That is, they randomly move in a cardinal direction, or remain still, with uniform probability.

\begin{figure}[htp!]
	\includegraphics[width=\textwidth]{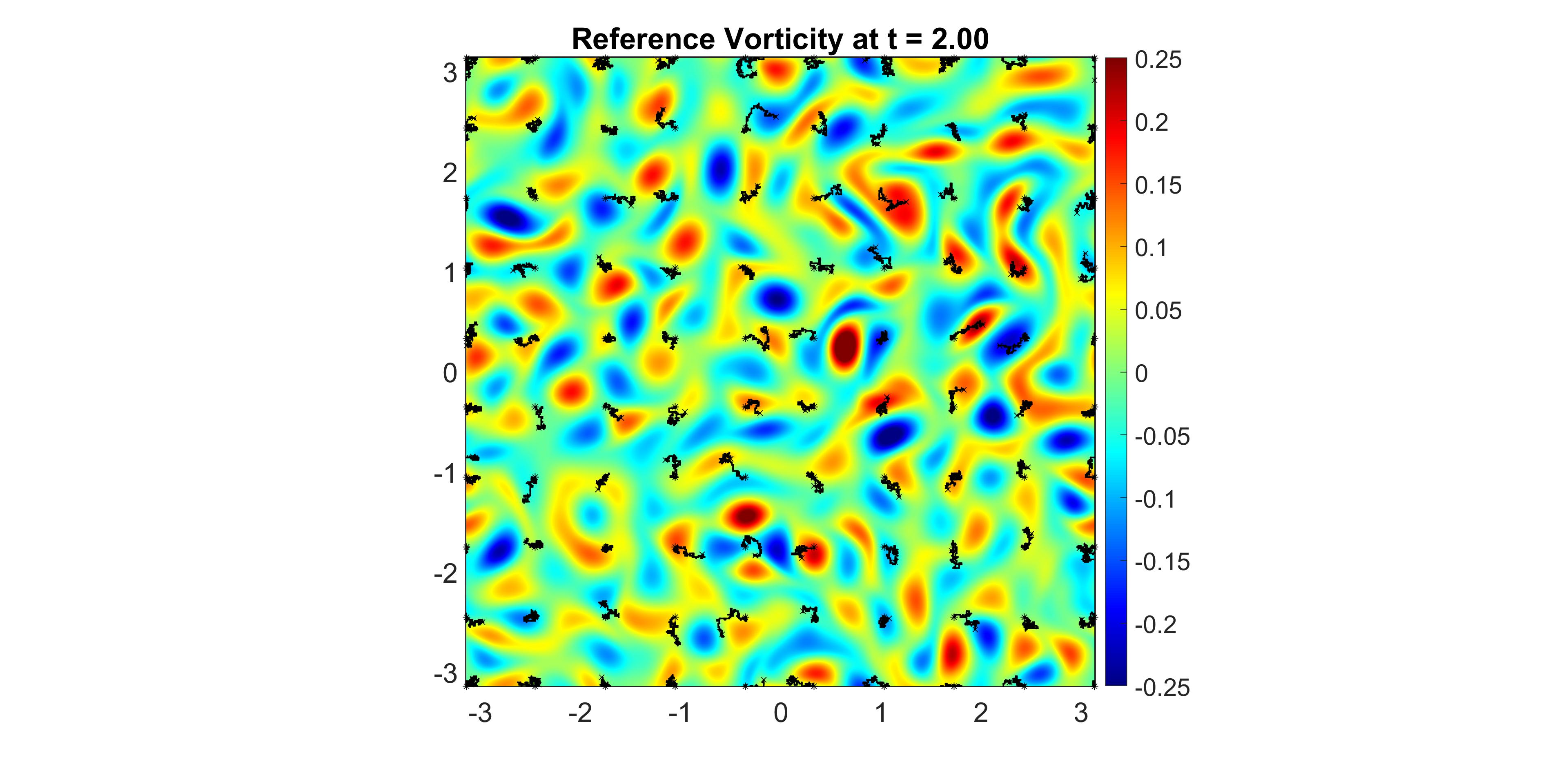}
	\caption{Graph displaying movement of observers using the creeps movement style displayed over the vorticity of the reference solution at time $t=2$. Observers were initially placed at 100 locations indicated by * and travel along a random walk indicated by black paths with ending locations marked by $\times$.}
	\label{fig:MovementCreep}
\end{figure}

\subsection{Lagrangian Particles}
The final movement strategy we considered is given by observers following Lagrangian particles trajectories. Observers are initialized into a uniform grid of spatial resolution $h>0$ and then follow the Lagrangian trajectory $\ell(t)$ as determined by 
\begin{align}
    \frac{d\ell }{dt} = u(\ell(t),t).
\end{align}
This is a very natural type of dynamic observer, as the movement of the observer comes exactly from the behavior of the system that we are attempting the capture in simulations of $v$. Once can consider sensors attached to buoys (e.g., Argo floats) set adrift in the ocean.

Lagrangian trajectories for observers were evolved forward using the third-order Adams-Bashforth method with values of $u$ at locations $(\ell(t),t)$. As in the case for random-sweeps, observers in the case of Lagrangian particles the observations do not in general occur at location on the underlying spatial grid.
Therefore we must approximate the value $u(\ell(t),t)$ via a linear interpolation of $u$. In this case we can consider $I_{h,t}(u-v)$ as follows:
\begin{align*}
    I_{h,t}(u-v) = I_{\ell(t),t}\circ J_{\ell(t)}(u-v)\\
\end{align*}
Here $J_{\ell(t)}(u-v)$ is a piecewise-linear interpolant used to approximate the values of $u-v$ in order to find the observational data and $I_{\ell(t),t}$ is another piecewise-linear interpolant used to extend the observational data to the numerical grid. 

\begin{figure}[htp!]
	\includegraphics[width=\textwidth]{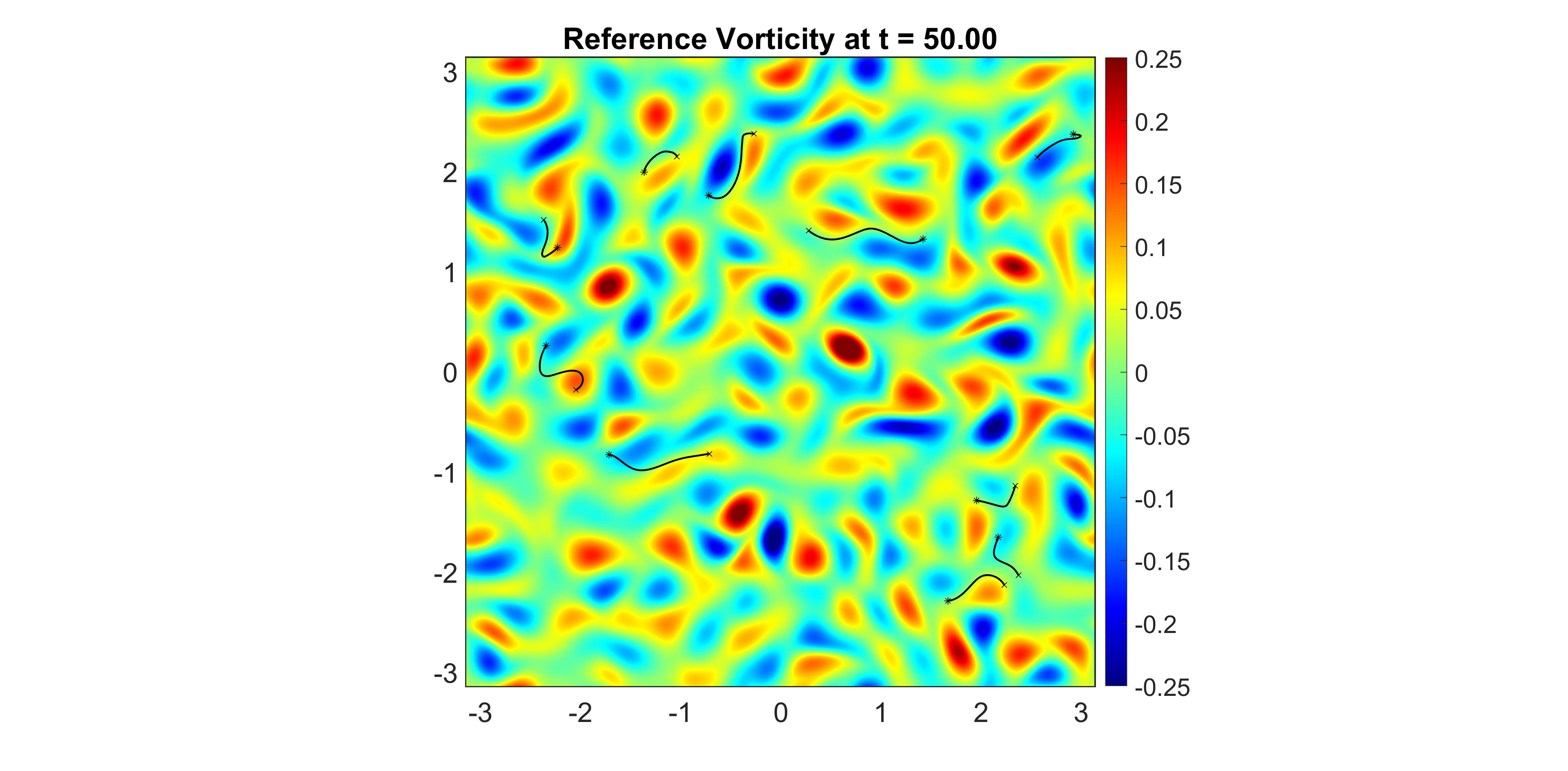}
	\caption{Graph displaying movement of observers given by Lagrangian particles displayed over the vorticity of the reference solution at time $t=50$. Observers were initially placed at 10 locations indicated by * and travel along random trajectories indicated by black paths with ending locations marked by $\times$.}
	\label{fig:MovementLagrange}
\end{figure}

\section{Computational Results}

All of our computations were done on the 2D NSE in the vorticity-stream function formulation using a fully dealiased pseudospectral code with physical domain $\nT^2 = [-\pi,\pi]^2$. 
As in \cite{Biswas_Bradshaw_Jolly2020} we used viscosity $\nu = 10^{-4}$  with forcing as given in \cite{Gesho_Olson_Titi_2015} that was multiplied by a scalar constant to obtain a forcing $f$ with Grashof number $G = 10^6$. Additionally, consistent with \cite{Biswas_Bradshaw_Jolly2020,Olson_Titi_2008_TCFD,Olson_Titi_2003}, we used a third-order Adams-Bashforth method with an integrating factor to solve for the linear term exactly. We used step size $\Delta t = 0.005$, and our initial data was generated by evolving zero initial data out to time $t = 25,000$. The spectrum of the initial data for the reference solution $u$ is plotted in Figure \Cref{fig:Spectrum} and one can see that spectrum is well-resolved (it also remained well-resolved for all times).  Note that the spectrum of the assimilation solution $v$ was temporarily slightly under-resolved in . Additionally, we used spatial resolution $N^2 = 1,024^2$ in our simulations in order to fully resolve the energy spectrum to machine precision (roughly $2.22\times10^{-16}$) before the 2/3's dealiasing cutoff (see, e.g., Figure \ref{fig:Spectrum}). 

Note that while the equation was simulated using the vorticity-stream function formulation, the interpolation term was calculated at the velocity level using nodal interpolation in physical space. To preserve the periodicity in the interpolation across the periodic ``boundary'' we interpolated using observed values over $[-3\pi, 3\pi]\times[-3\pi,3\pi]$ instead of just over $\nT^2 = [-\pi,\pi]\times[-\pi,\pi]$, and then truncated the result back to the physical domain. 
Since we are simulating the equations at the vorticity-stream function level, we calculate  $\|\psi - \tilde{\psi}\|_{L^2}$ as the error between the simulated and reference solutions. We additionally calculate both $\norm{\omega - \tilde{\omega}}_{L^\infty}$ and $\norm{\omega - \tilde{\omega}}_{L^2}$ (see   \cref{fig:SIM_min_psi,fig:CPU_min_psi,fig:SIM_min_omega1,fig:CPU_min_omega1,fig:SIM_min_omega2,fig:CPU_min_omega2,fig:SIM_equal_psi,fig:CPU_equal_psi,fig:SIM_equal_omega1,fig:CPU_equal_omega1,fig:SIM_equal_omega2,fig:CPU_equal_omega2}). Here $\psi$ denotes the stream function of the $u$ equation with $\omega$ being the vorticity of $u$, and  $\tilde{\cdot}$ denoting the corresponding quantity in system \eqref{NSE-v}.  

Using the numerical methods discussed above, we simulated dynamic observers that move in accordance with the movement styles. 
We found that the choice of $\mu$ varied for different movement methods.
For the uniform static grid it was found that $\mu = 10$ was the optimal choice of $\mu$. It appears that the choice of larger $\mu$ values is restricted by a CFL condition on the time-step; namely $\Delta t<2/\mu$. We observed, similar to our previous work, \cite{Larios_Victor_2019} that highly mobile movement schemes, such as thin-sweeps, do not appear to be subject to the same restriction on $\mu$, perhaps due to the numerically destabilizing effect of large $\mu$-values being only briefly concentrated in any given spatial location.  For simulations of sweeps, we used the value $\mu = 30$, which was unstable in the case of the static uniform grid of observers. In our simulations we observed that larger values of $\mu$ tended to lead to faster convergence in the case of thin-sweeps. In simulations, we tested thick-sweeps and found that the method did allow for larger choices of $\mu$ when the speed of the observer movement was increased to that of the thin-sweeps (from $1\frac{\Delta x}{\Delta t}$ to $3\frac{\Delta x}{\Delta t}$), further corroborating the ``brief local destabiliztion'' explanation discussed above. While thick-sweeps did allow for a larger choice of $\mu$ we did not observe faster convergence rates in these trials. We conjecture that, at least typically, the larger values of $\mu$ are optimal for very fast moving observer regimes in which any given point will only be observed for a short amount of time.

In preliminary simulations, we observed that certain $\mu$ values for which the static uniform grid remained stable could introduce instabilities in mobile observer schemes (except in the case of thin-sweeps and thick-sweeps). 
This is perhaps due to a CFL condition restricting the length-scale $h$, the minimum length between observers. Since $h$ is not restricted for most of the mobile observer paradigms, $\mu$ values should be smaller to avoid this instability. 
We found that the value $\mu = 10$ was optimal for the uniform static grid and did not cause instability in the mobile observers for our choice of spatial resolution $N^2=1024^2$. In preliminary simulations with spatial resolution $N^2 = 512^2$ we saw that the mobile observer methods required $\mu$ values smaller than the the optimal choice of $\mu$ for the uniform static grid to avoid instability. However this instability was not seen in any of the higher resolution tests. 
While $\mu$ can potentially be further optimized, as generally higher $\mu$ values coincide with faster decay in the error, it is worth noting that our results showed faster rates of convergence using the same $\mu$ values, but with smaller numbers of observers required, (see \cref{fig:SIM_min_psi,fig:CPU_min_psi}).

To test these methods computationally, we first ran simulations for the uniform static grid case to establish a baseline for comparison. We found that, for our choices of parameters, a uniform grid of $5,625$ observers would converge to machine precision $10^{-14}$ by approximately time $t = 70$ (see \cref{fig:SIM_min_psi,fig:CPU_min_psi,fig:SIM_min_omega1,fig:CPU_min_omega1,fig:SIM_min_omega2,fig:CPU_min_omega2,fig:SIM_equal_psi,fig:CPU_equal_psi,fig:SIM_equal_omega1,fig:CPU_equal_omega1,fig:SIM_equal_omega2,fig:CPU_equal_omega2}). Using this, we varied the number of observers in the case of the bleeps, the sweeps, and the creeps in order to obtain convergence at approximately the same time ($t = 70$) for the sake of comparison. The results from these studies can be seen in \cref{fig:SIM_min_psi,fig:CPU_min_psi,fig:SIM_min_omega1,fig:CPU_min_omega1,fig:SIM_min_omega2,fig:CPU_min_omega2,fig:SIM_equal_psi,fig:CPU_equal_psi,fig:SIM_equal_omega1,fig:CPU_equal_omega1,fig:SIM_equal_omega2,fig:CPU_equal_omega2}. As we see from these results, all three of these methods enjoy  exponential convergence to the reference solution in the $L^2$ norm. Moreover, all three methods converge to machine precision faster than the uniform grid using significantly fewer observers. The bleeps in particular demonstrated almost an order of magnitude improvement in the number of observers compared to the uniform grid. 

Next, we ran simulations utilizing mobile observers with each initialized to approximately 5,625 observers to compare with the static uniform grid of 5,625 observers. As expected, we observed that all of the mobile observer methods converged exponentially fast to the reference solution $u$ at a faster rate than the uniform static grid. Our results can be seen in \cref{fig:SIM_min_psi,fig:CPU_min_psi,fig:SIM_min_omega1,fig:CPU_min_omega1,fig:SIM_min_omega2,fig:CPU_min_omega2,fig:SIM_equal_psi,fig:CPU_equal_psi,fig:SIM_equal_omega1,fig:CPU_equal_omega1,fig:SIM_equal_omega2,fig:CPU_equal_omega2}.


One additional metric we used for comparison with the uniform static grid was the CPU time required to obtain convergence to machine precision. It is important to check that the extra complexity built into simulating moving observers does not make these methods prohibitively expensive to simulate when compared to the static uniform grid. In all of the simulations described above we also captured the CPU time. Graphs of the error vs the CPU time can be seen in \cref{fig:SIM_min_psi,fig:CPU_min_psi,fig:SIM_min_omega1,fig:CPU_min_omega1,fig:SIM_min_omega2,fig:CPU_min_omega2,fig:SIM_equal_psi,fig:CPU_equal_psi,fig:SIM_equal_omega1,fig:CPU_equal_omega1,fig:SIM_equal_omega2,fig:CPU_equal_omega2}. Note that the CPU time required to update the locations is included in these calculations of total CPU time. This includes the time to randomly determine locations, move the rectangular subdomain, and determine random walks for bleeps, sweeps, and creeps, respectively (this would likely not be needed when working with real-world data, but on the other hand, working real-world data may bring other CPU-intensive operations, such as I/O costs or off-grid interpolation on non-uniform grids).  Additionally, note that since the thin-sweeps strategy use the exact solution over a small rectangular subdomain there is no interpolation performed. This leads to a very low CPU time in those figures. In the plots, we found that when the number of observers is chosen to ensure convergence by (roughly) time $t = 70$, all of the dynamic observer movement strategies except Lagrangian particles obtain convergence significantly faster in CPU time than the uniform grid case. This is perhaps due to the fewer number of observers used,  simplifying the calculations of the interpolation term. Additionally, when all of the methods are simulated using approximately $5,625$ observers, the mobile observer methods still obtain faster convergence than the static uniform grid. In this case, it is likely due to the rate of convergence, which can be seen in \cref{fig:SIM_min_psi,fig:CPU_min_psi,fig:SIM_min_omega1,fig:CPU_min_omega1,fig:SIM_min_omega2,fig:CPU_min_omega2,fig:SIM_equal_psi,fig:CPU_equal_psi,fig:SIM_equal_omega1,fig:CPU_equal_omega1,fig:SIM_equal_omega2,fig:CPU_equal_omega2}.

We note that all of the data assimilation methods, except thin and thick sweeps, show a steep decay in the error initially until time $t\approx1$. This steep decay is possibly due to the global distribution of observers across the domain. This steeper decay can be seen in the case of thick-sweeps when initialized with 1,400 observers moving with constant velocity $3\frac{\Delta x}{\Delta t}$ in the $x$ direction until approximately time $t = 8$ (see \cref{fig:SIM_min_psi}). Moreover we note that this steep decay was not seen for thick-sweeps with the constant velocity $1\frac{\Delta x}{\Delta t}$ in the $x$ direction (see \cref{fig:SIM_thick_comparison}).

In the above simulations it appears that performing thin-sweeps is a good strategy for mobile data assimilation however it requires the exact solution in a small rectangle. All of the other movement schemes are adaptable to higher resolution with no difficulties, but with thin-sweeps for a fixed rectangular region the number of observers scales with the resolution.
Physically this may be feasible if the data collection comes from a device that gathers a velocity profile along one dimension, however it is not feasible if the velocity is measured using observers placed at discrete points. This lead to the generalization of thin-sweeps to thick-sweeps, which allows for mobile observers that sample locally in space and move continuously in time without requiring a number of observers proportional to the spatial resolution of the underlying numerical grid. 
Here, instead of having access to the full solution in a rectangle, we extend the rectangular area to cover a quarter of the domain and interpolate locally over a uniform grid of observers distributed locally on this subdomain. In our simulations we found that the velocity of the observers was an important factor in the convergence rate (see \Cref{fig:SIM_thick_comparison,fig:CPU_thick_comparison}). One may be able to optimize the velocity of the observers based on the size of the domain being sampled to improve convergence rates further, but we do not explore this in the present work.

\section{Conclusions}
Our simulations of mobile observers provide evidence that highly mobile observers such as the bleeps give the fastest convergence rates using the fewest number of observers when compared to a static uniform grid of observers. All mobile methods demonstrated convergence faster in simulation time than the static uniform grid when initialized with the same number of observers; however, the methods that demonstrated the highest convergence rates were the bleeps and random-sweeps, both of which require some level of randomness.
While sensors with random movement may be difficult to implement physically, this does demonstrate the need for observing data throughout the domain. 



While all of the methods converged to machine precision faster in simulation time than the static uniform grid, the thick-sweeps and Lagrangian mobile observer movement schemes performed poorly in comparison to the rest of the mobile methods.
In particular, neither thick-sweeps nor Lagrangian particle observers demonstrate the improvements in CPU time seen in the other mobile methods. 
In the case of thick-sweeps, the poor performance relative to the other mobile methods is likely due to lack of mobility. While thick-sweeps do move observers continuously, all of the simulations done in this study used $b = 1$ which results in slow movement across the domain. Larger $b$ values would likely make thick-sweeps behave similarly to thin-sweep in convergence rates. Similarly, we believe that the Lagrangian particle observers perform poorly due to the slow movement of fluid particles and the tendency for particles to become trapped in vorticity structures (e.g., eddies).

An additional point of note is the performance of thin-sweeps. 
Recall that thin-sweeps measures the full solution within a small rectangle and inserts this directly into the $v$ equation without interpolation of any sort. This leads to a drastic decrease in CPU time when compared to the other methods which can be seen in \cref{fig:CPU_min_psi,fig:CPU_min_omega1,fig:CPU_min_omega2,fig:CPU_equal_psi,fig:CPU_equal_omega1,fig:CPU_equal_omega2,fig:SIM_min_psi,fig:SIM_min_omega1,fig:SIM_min_omega2,fig:SIM_equal_psi,fig:SIM_equal_omega1,fig:SIM_equal_omega2}. Additionally, due to the requirement of measuring the full solution inside a small rectangle, the number of observers is misleading. As the resolution increases, the thin-sweep strategy will need an increasing amount of observers if the rectangular area remains the same size. The need to measure the full solution essentially means that the thin-sweep strategy requires an infinite number of observers. As currently implemented, the number of observers determines the width of the rectangular area, which for $3072$ is $3072/N \Delta x = 3\Delta x$. Another restriction of thin-sweeps is that this requirement that the number of observers must be a multiple of the resolution $N$. 

Due to the shortcomings of both thin-sweeps and thick-sweeps, we propose that a hybrid of the two methods would be practical to implement and would have some of the beneficial properties of both while being more readily implementable physically. This hybrid method would replace the rectangular area of thick-sweeps with the much smaller rectangle from thin-sweeps. Additionally, as in thin-sweeps, this hybrid method would keep the movement of the rectangular area proportional to its width. However, unlike thin-sweeps, we would actually interpolate over the area within the rectangle instead of requiring the full solution. We will explore these ideas in a future work.

From our simulations it appears that the bleeps method requires the smallest number of observers to achieve convergence comparable to the uniform static grid. In the description of the bleeps movement strategy, we noted that the locations of observers are changed at each timestep. In our previous simulations, it was found that the bleep strategy also demonstrates exponential convergence to the reference solution when observers are placed in random locations that changed after a specified amount of time. We noted that the error decreased sharply for a short amount of time after the observers were moved. In order to maintain this steep rate of decay, we chose to randomize the observer locations at every timestep. This suggests that one could optimize the randomization time, which will be the subject of a future work.

When considering the viability of mobile methods, it is important to consider the feasibility of physical implementations. In particular, one should ask if one can  realize these methods using actual sensors physically moving in time. Moreover, in physical implementations, mobile observers may be more expensive than non-mobile observers and more prone to break down and potentially suffer measurement error.  
Therefore, factoring in cost of number of observers should be taken into account in real-world applications, but this is beyond the scope of the present work.
In this study, we assumed that all observers had the same level of accuracy (namely, they are assumed to have perfect accuracy).  Although studies of the AOT algorithm have been carried out in the context of noisy data \cite{Bessaih_Olson_Titi_2015,Carlson_Hudson_Larios_Martinez_Ng_Whitehead_2021}, showing that the algorithm with static observations is robust with respect noisy measurements, factoring in different levels of accuracy for different observers may be crucial in real-world implementations.  
(For example, sensors on airplanes tend to be much more accurate, but more expensive that ground-based sensors.)

We note that, although we made an effort to follow good coding practices, no specific effort was made to optimize the computation of the various moving observer schemes, nor their interpolation, nor the nudging constants $\mu$, except in the case of static observers.  However, the moving-observer schemes still outperform the static observer case in simulation time.  When measured in CPU time the moving-observer schemes outperform the static observer case except in the case of Langrangian observers.  It appears that Lagrangian tracer particles tend to get stuck in eddies, decreasing their global coverage.  However, random the random sampling locations of the bleeps demonstrated significantly better error.  Therefore, for example, in the case of observational devices attached to buoys, it may be worth considering mechanisms to ``kick'' the buoys into a different location.  For example, a small motor or wind sail attached to Argo buoys could lead to significant decreases the error.

In summary, this work demonstrates the enormous potential for improving data assimilation speed, accuracy, and the number of observers required by allowing for observer to move through the domain in certain patterns, sampling wider regions of the flow.  These results open the door to many new possibilities, including optimizing movement patterns (using, e.g., machine learning strategies), combining data from moving and non-moving sources (using, e.g., data from ground-based stationary sensors and fast-moving drones), heterogeneous data sources (e.g., temperature and salinity of the ocean), and many others.  We will explore these ideas in future works.

\appendix
\section{Additional Figures}

\begin{figure}[htp!]
	\includegraphics[width=\textwidth]{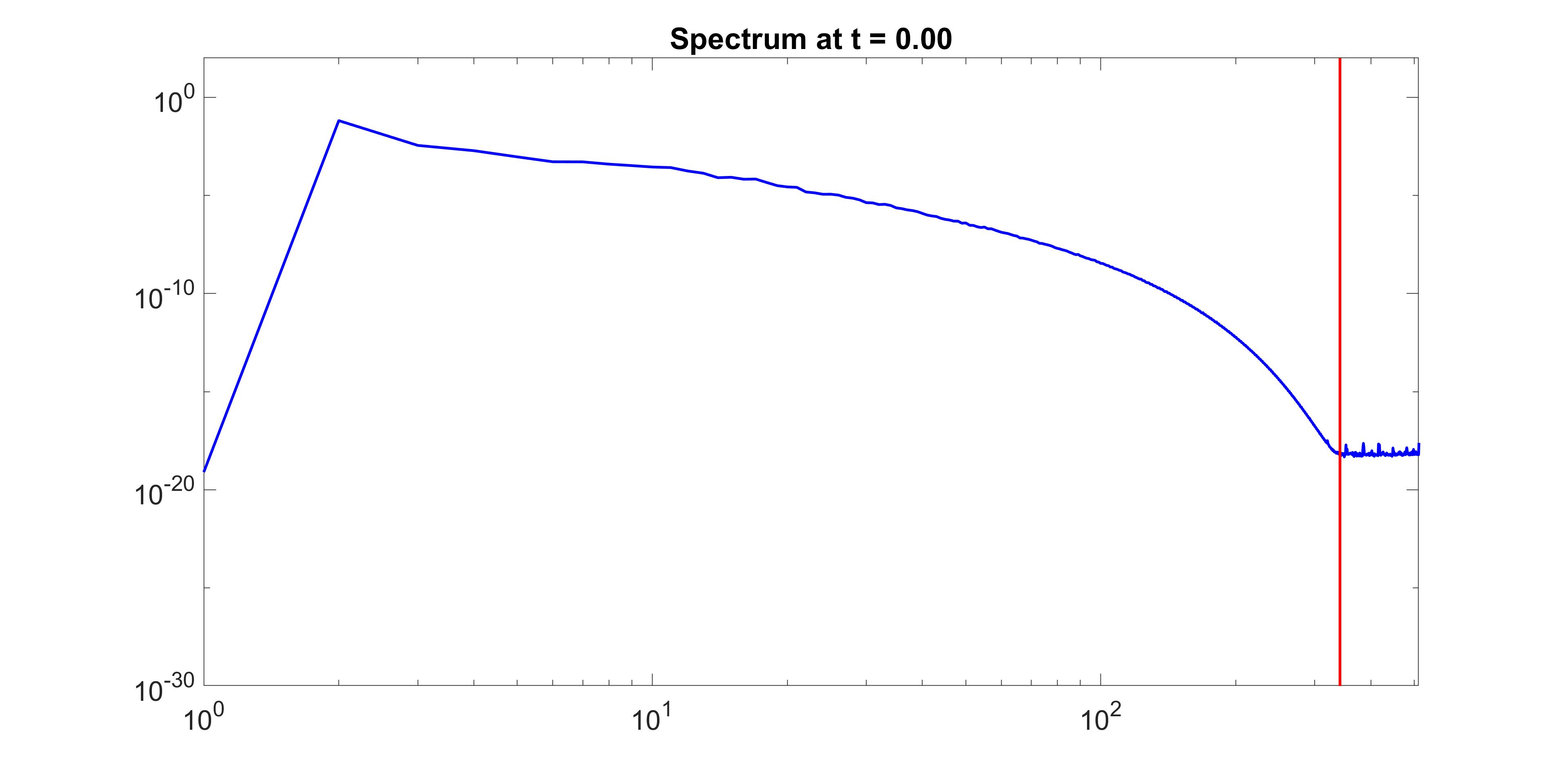}
	\caption{Energy spectrum of the initial data with $\nu = 10^{-4}$, $G = 10^6$, and $\Delta t = 0.005$. The vertical red line is the 2/3's dealiasing cutoff at $\frac23\frac{N}{2}=341.\overline{3}$.}
	\label{fig:Spectrum}
\end{figure}

\begin{figure}[htp!]
	\includegraphics[width=.8\textwidth]{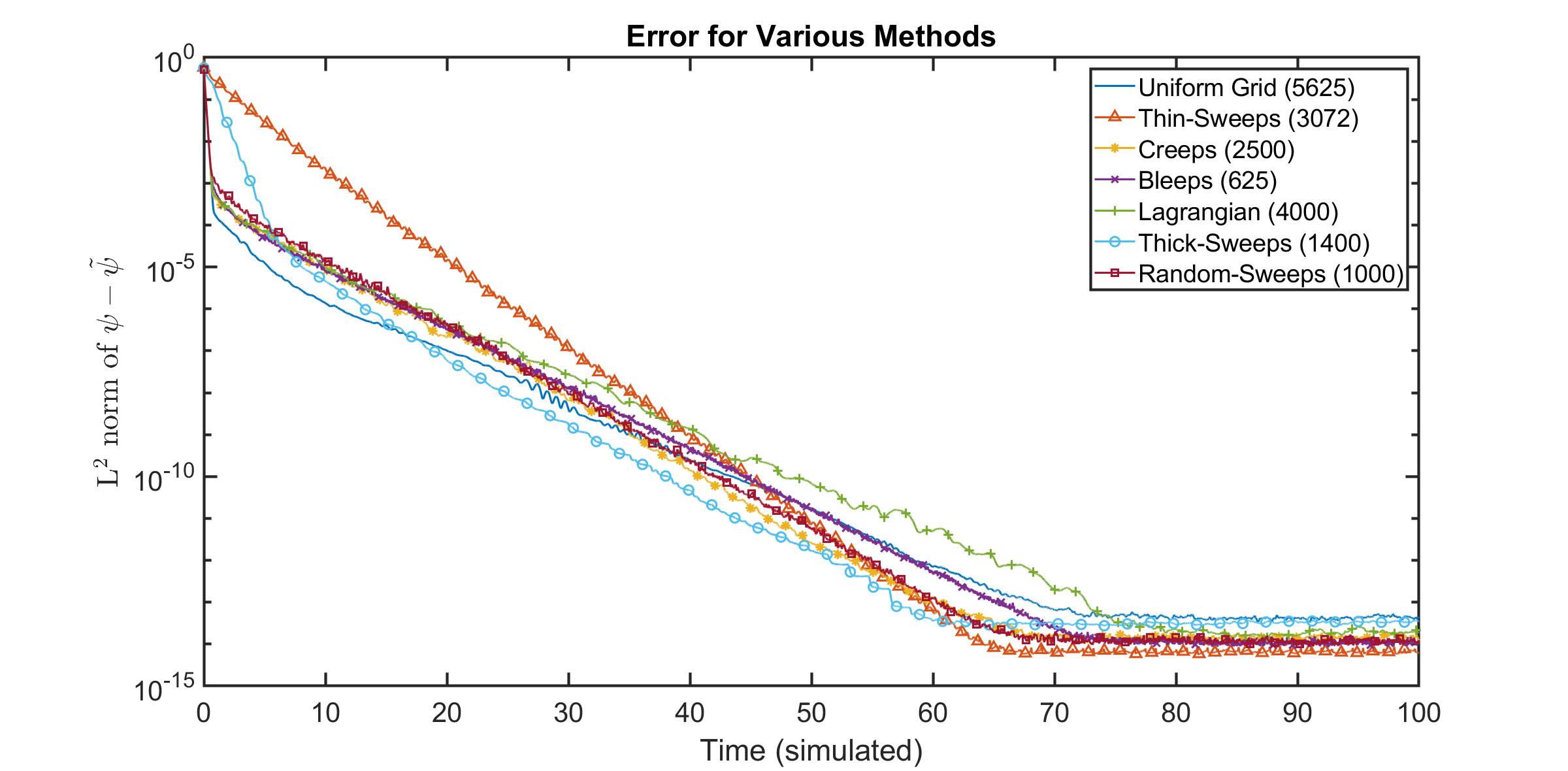}
	\caption{Comparison of error over time for static uniform grid and dynamic observer strategies with number of observers optimized to ensure convergence around the time $t = 70$ (log-linear plot). 
	See legend entries for the number of observers used in each method. $\mu = 10$ for all methods except thin-sweeps, which uses $\mu = 30$. Thick-sweeps has constant $x$ velocity $3\frac{\Delta x}{\Delta t}$.
	}
	\label{fig:SIM_min_psi}
\end{figure}

\begin{figure}
	\includegraphics[width=\textwidth]{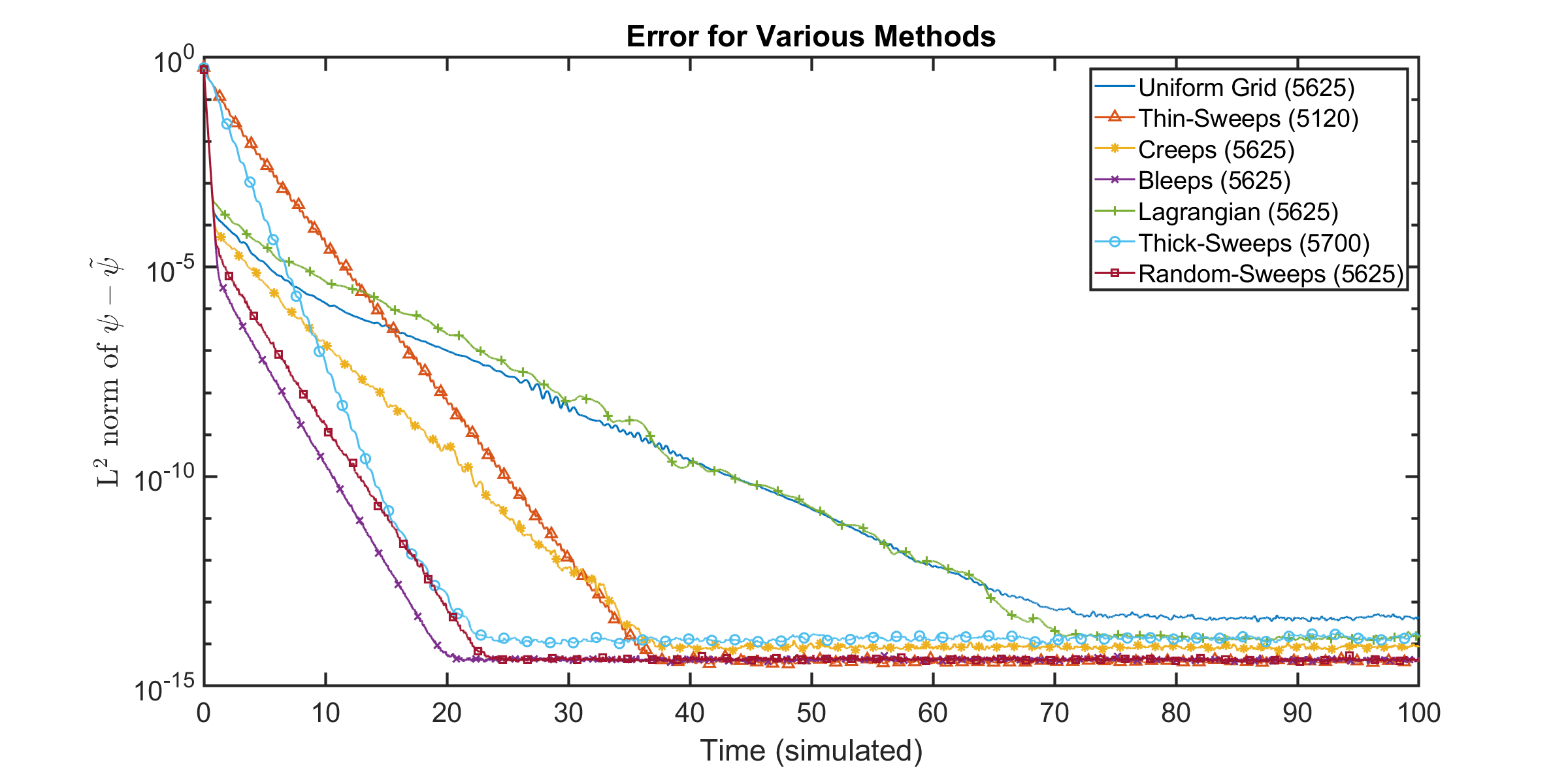}
	\caption{Comparison of error over time for static uniform grid and dynamic observer strategies all initialized using approximately 5,625 observers (log-linear plot). 
	See legend entries for the number of observers used in each method. $\mu = 10$ for all methods except thin-sweeps, which uses $\mu = 30$. Thick-sweeps has constant $x$ velocity $3\frac{\Delta x}{\Delta t}$.
	}
	\label{fig:SIM_equal_psi}
\end{figure}

\begin{figure}[htp!]
	\includegraphics[width=\textwidth]{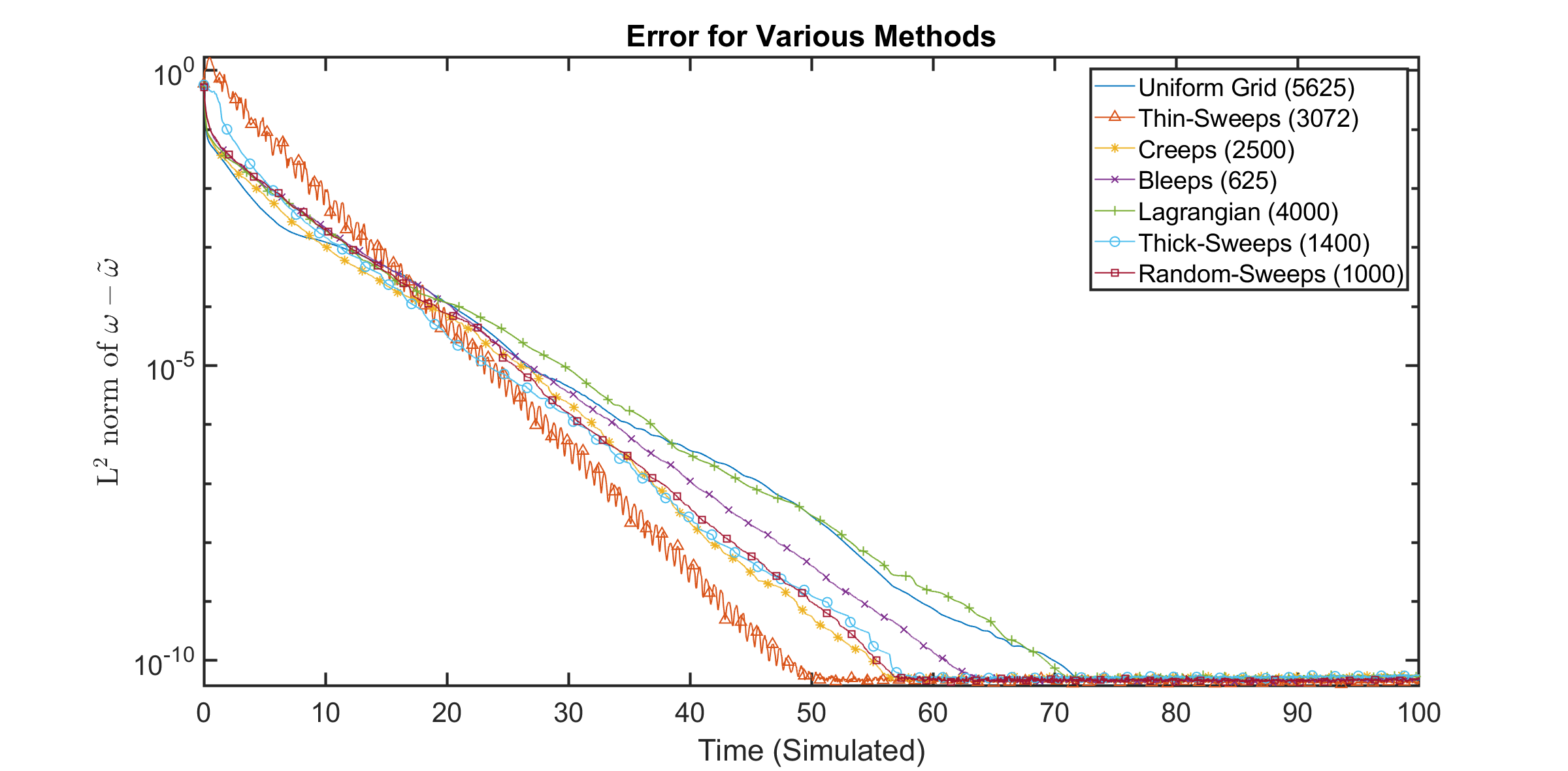}
	\caption{Comparison of error over time for static uniform grid and dynamic observer strategies with number of observers optimized to ensure convergence around the time $t = 70$ (log-linear plot). 
	See legend entries for the number of observers used in each method. $\mu = 10$ for all methods except thin-sweeps, which uses $\mu = 30$. Thick-sweeps has constant $x$ velocity $3\frac{\Delta x}{\Delta t}$.
	}
	\label{fig:SIM_min_omega1}
\end{figure}

\begin{figure}
	\includegraphics[width=\textwidth]{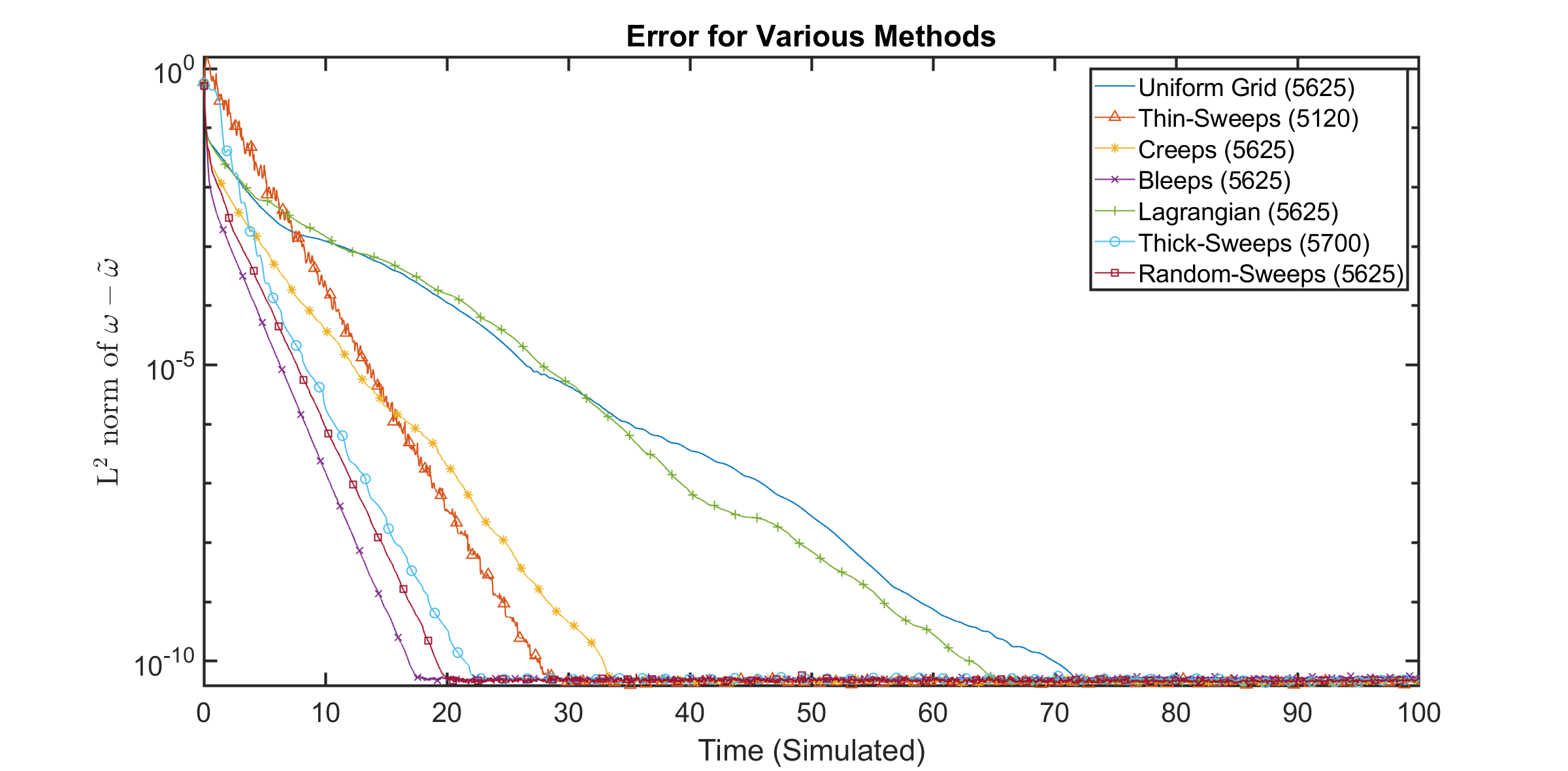}
	\caption{Comparison of error over time for static uniform grid and dynamic observer strategies all initialized using approximately 5,625 observers (log-linear plot). 
	See legend entries for the number of observers used in each method. $\mu = 10$ for all methods except thin-sweeps, which uses $\mu = 30$. Thick-sweeps has constant $x$ velocity $3\frac{\Delta x}{\Delta t}$.
	}
	\label{fig:SIM_equal_omega1}
\end{figure}

\begin{figure}[htp!]
	\includegraphics[width=\textwidth]{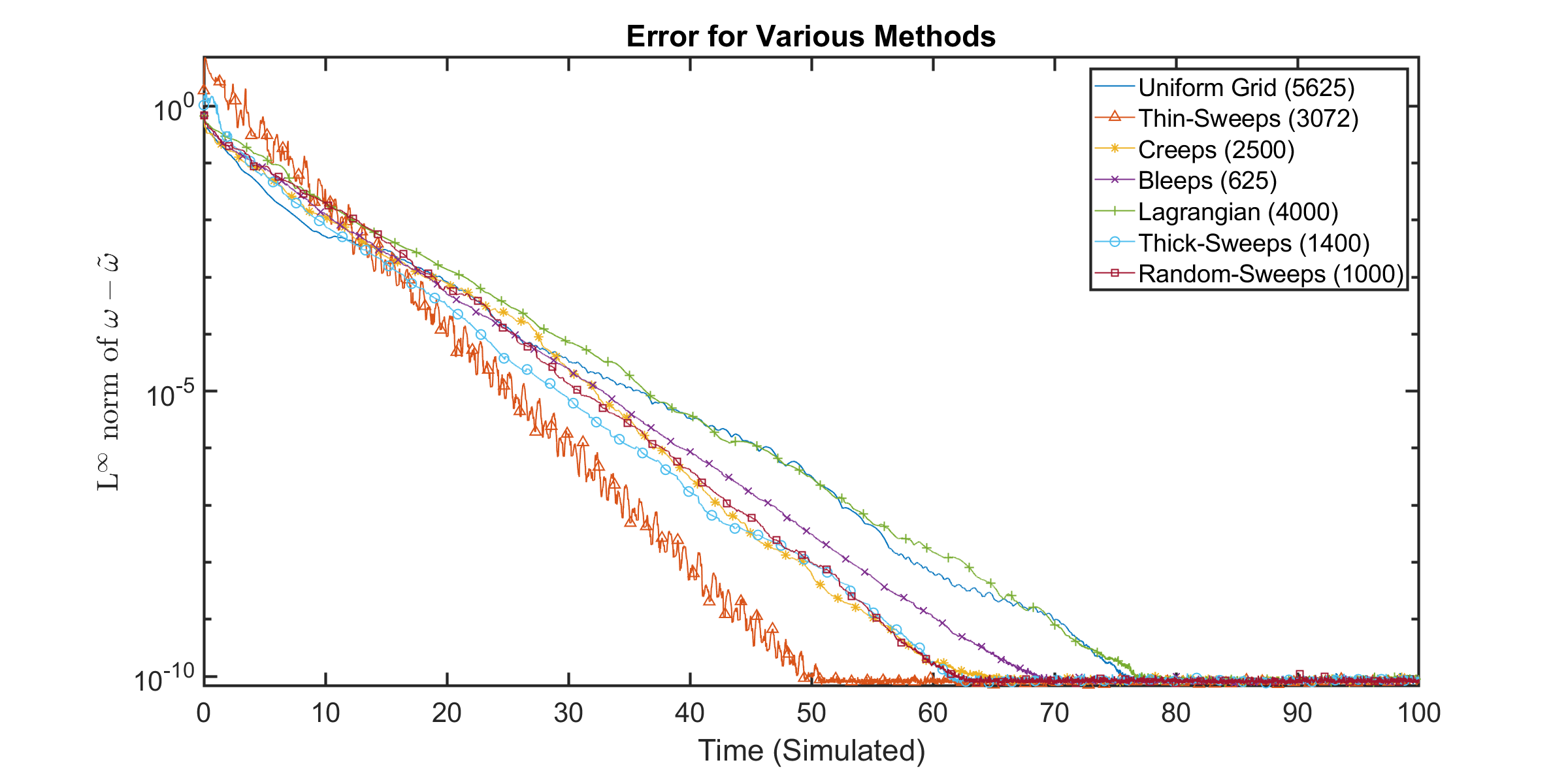}
	\caption{Comparison of error over time for static uniform grid and dynamic observer strategies with number of observers optimized to ensure convergence around the time $t = 70$ (log-linear plot).
	See legend entries for the number of observers used in each method. $\mu = 10$ for all methods except thin-sweeps, which uses $\mu = 30$. Thick-sweeps has constant $x$ velocity $3\frac{\Delta x}{\Delta t}$.
	}
	\label{fig:SIM_min_omega2}
\end{figure}

\begin{figure}
	\includegraphics[width=\textwidth]{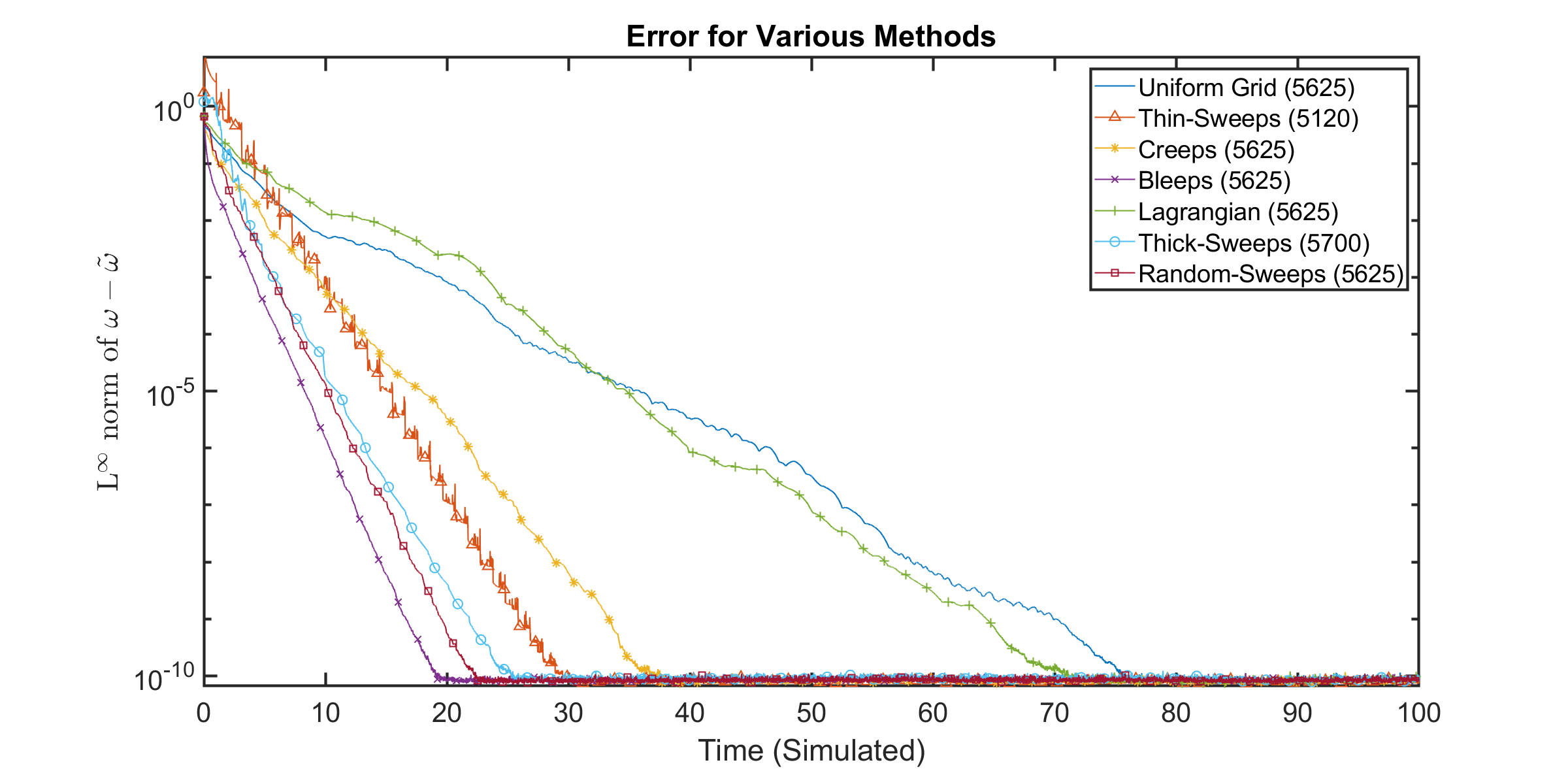}
	\caption{Comparison of error over time for static uniform grid and dynamic observer strategies all initialized using approximately 5,625 observers (log-linear plot). 
	See legend entries for the number of observers used in each method. $\mu = 10$ for all methods except thin-sweeps, which uses $\mu = 30$. Thick-sweeps has constant $x$ velocity $3\frac{\Delta x}{\Delta t}$.
	}
	\label{fig:SIM_equal_omega2}
\end{figure}

\begin{figure}
	\includegraphics[width=\textwidth]{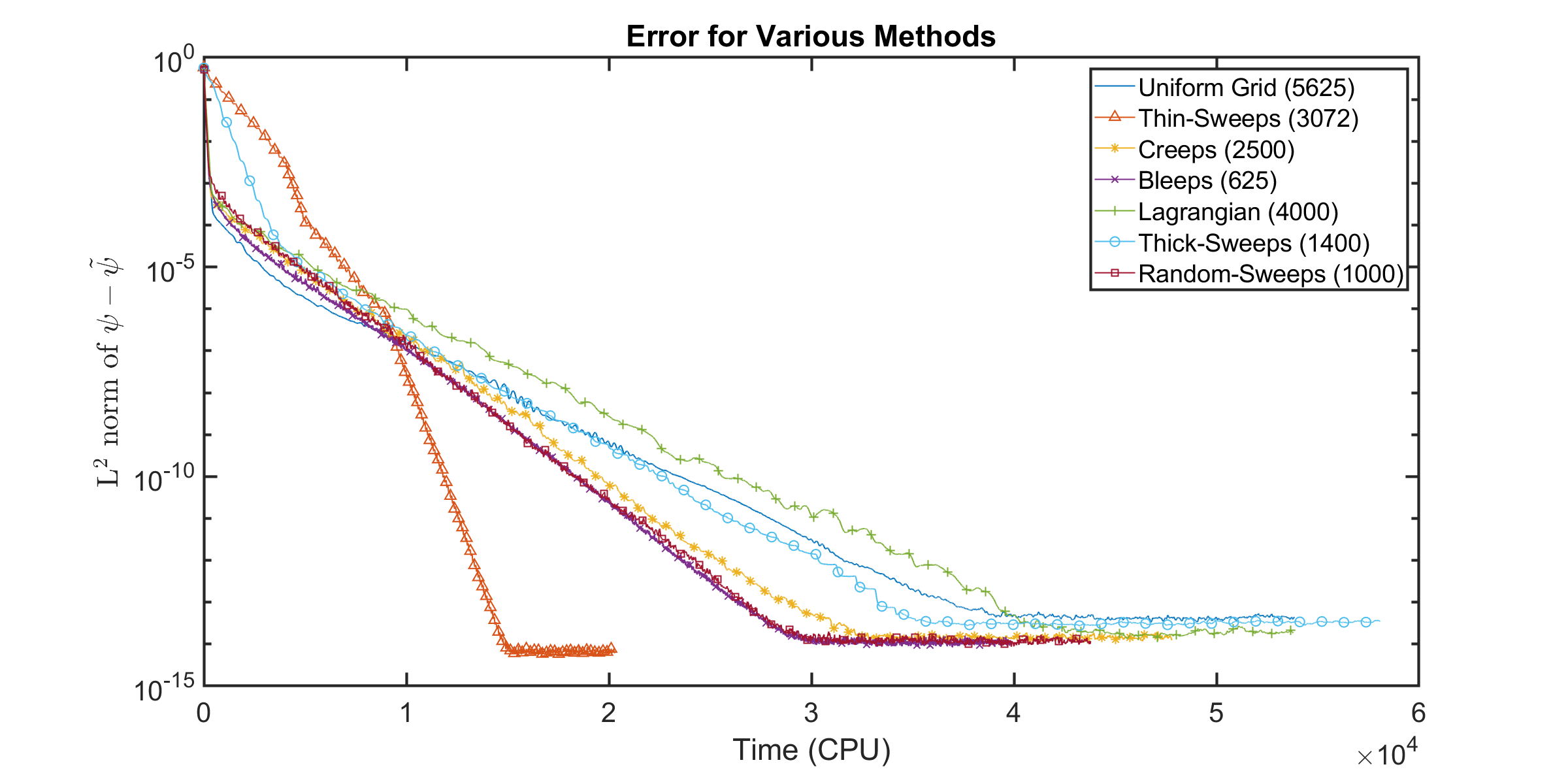}
	\caption{Comparison of error over CPU time for static uniform grid and dynamic observer strategies with number of observers optimized to ensure convergence around the time $t = 70$ (log-linear plot).
	See legend entries for the number of observers used in each method. $\mu = 10$ for all methods except thin-sweeps, which uses $\mu = 30$. Thick-sweeps has constant $x$ velocity $3\frac{\Delta x}{\Delta t}$.
	}
	\label{fig:CPU_min_psi}
\end{figure}

\begin{figure}
	\includegraphics[width=\textwidth]{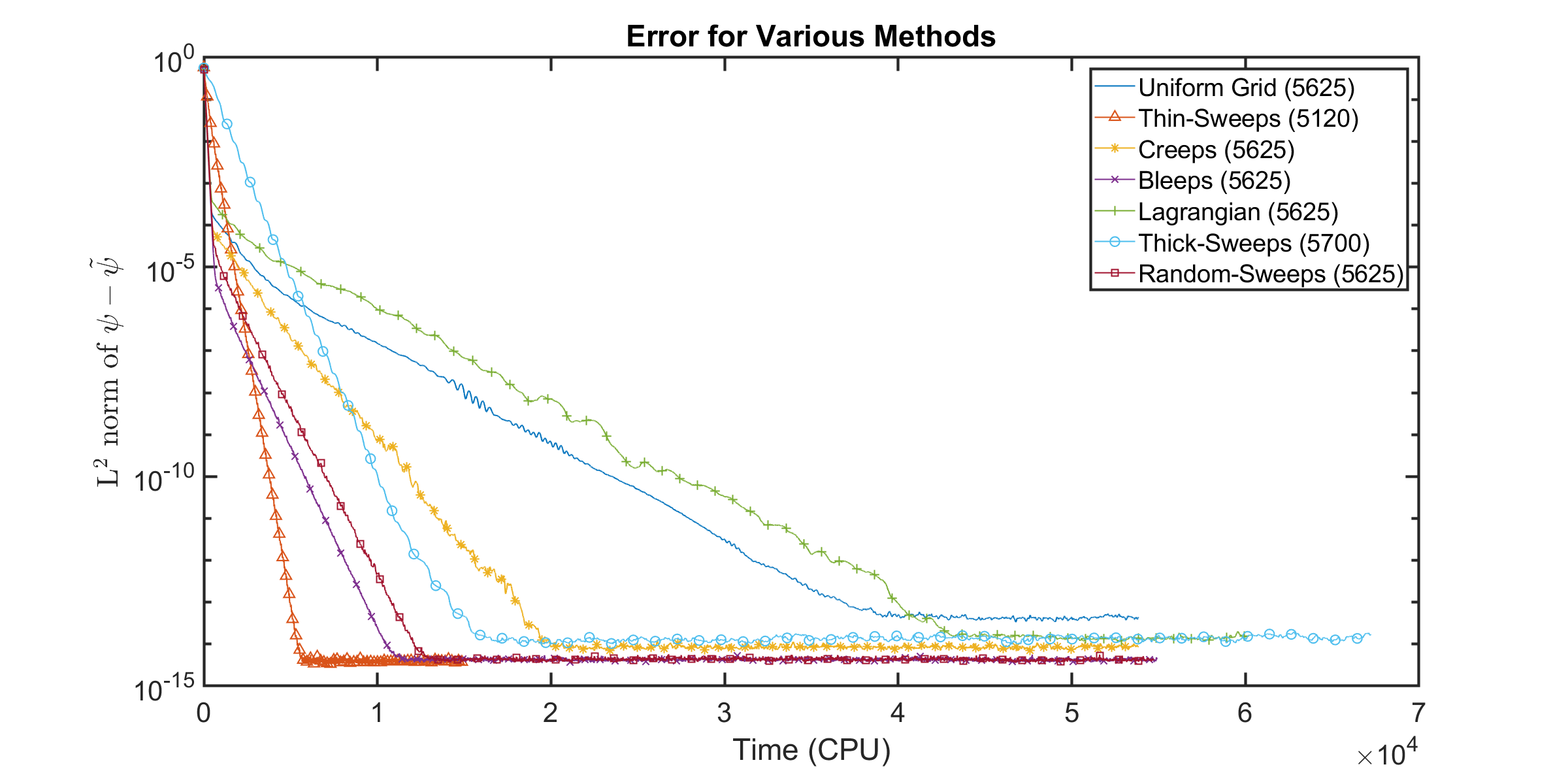}
	\caption{Comparison of error over CPU time for static uniform grid and dynamic observer strategies all initialized using approximately 5,625 observers (log-linear plot). 
	See legend entries for the number of observers used in each method. $\mu = 10$ for all methods except thin-sweeps, which uses $\mu = 30$. Thick-sweeps has constant $x$ velocity $3\frac{\Delta x}{\Delta t}$.
	}	\label{fig:CPU_equal_psi}
\end{figure}

\begin{figure}
	\includegraphics[width=\textwidth]{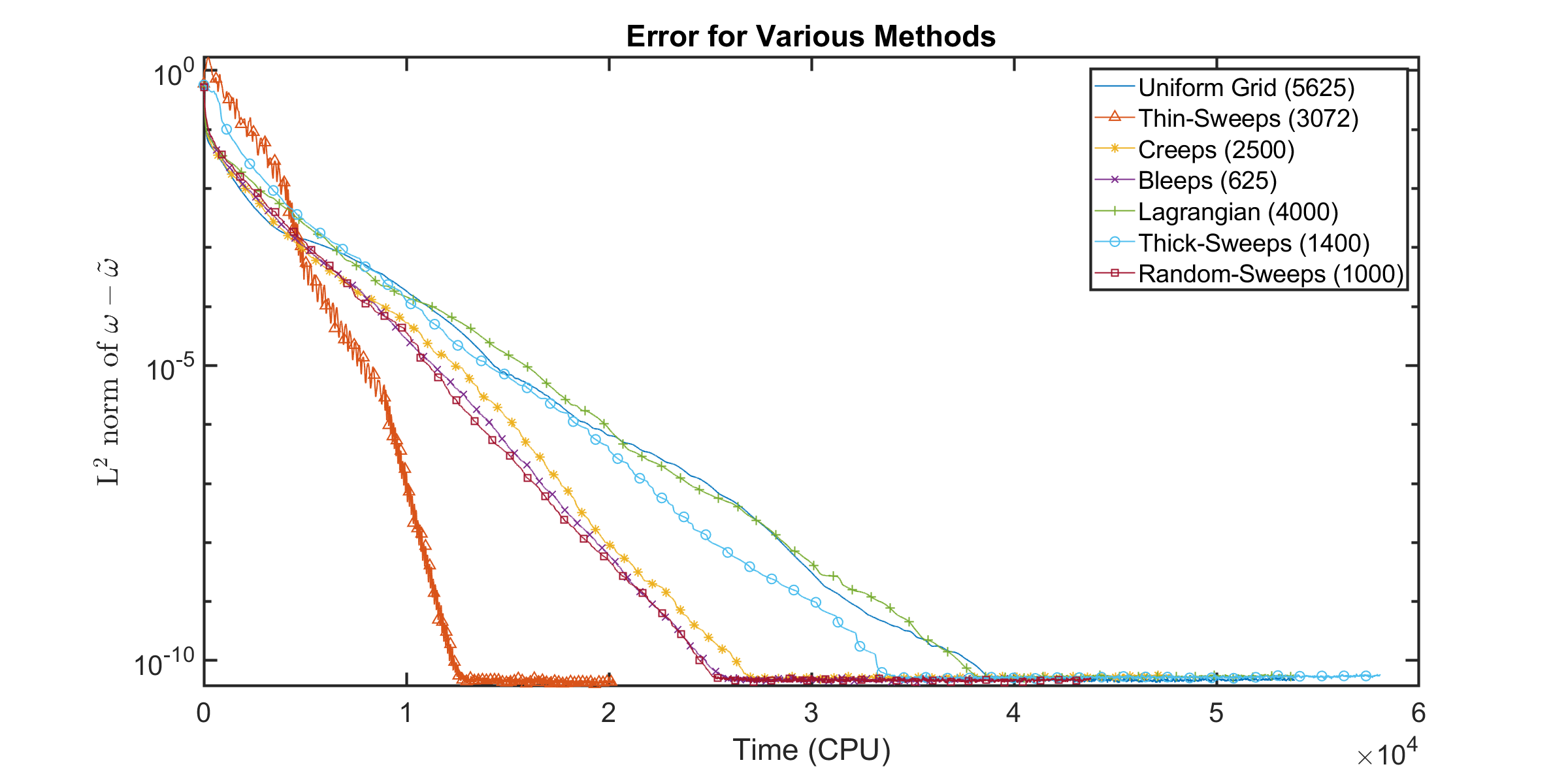}
	\caption{Comparison of error over CPU time for static uniform grid and dynamic observer strategies with number of observers optimized to ensure convergence around the time $t = 70$ (log-linear plot). 
	See legend entries for the number of observers used in each method. $\mu = 10$ for all methods except thin-sweeps, which uses $\mu = 30$. Thick-sweeps has constant $x$ velocity $3\frac{\Delta x}{\Delta t}$.
	}
	\label{fig:CPU_min_omega1}
\end{figure}

\begin{figure}
	\includegraphics[width=\textwidth]{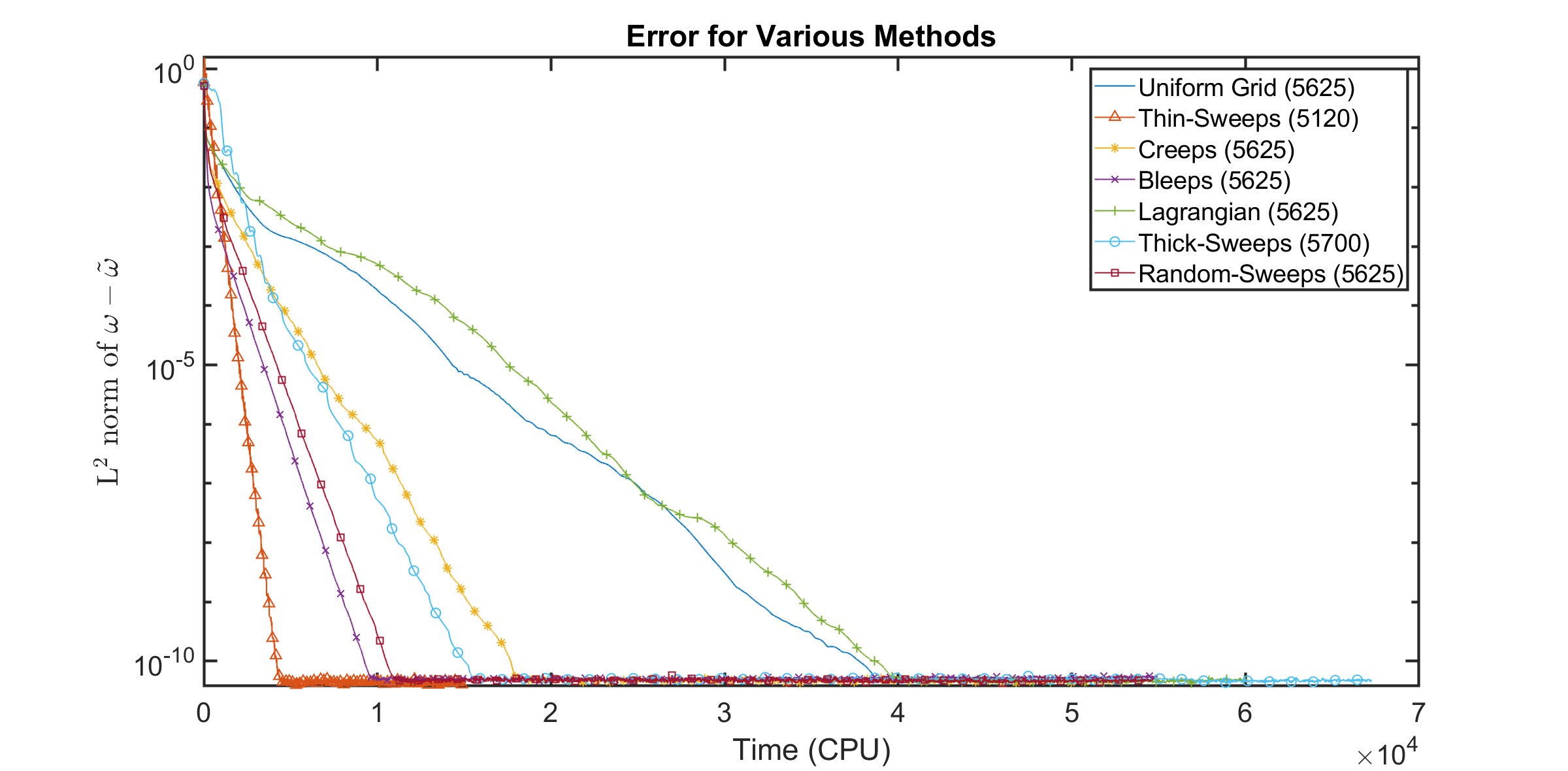}
	\caption{Comparison of error over CPU time for static uniform grid and dynamic observer strategies all initialized using approximately 5,625 observers (log-linear plot). 
	See legend entries for the number of observers used in each method. $\mu = 10$ for all methods except thin-sweeps, which uses $\mu = 30$. Thick-sweeps has constant $x$ velocity $3\frac{\Delta x}{\Delta t}$.
	}	\label{fig:CPU_equal_omega1}
\end{figure}

\begin{figure}
	\includegraphics[width=\textwidth]{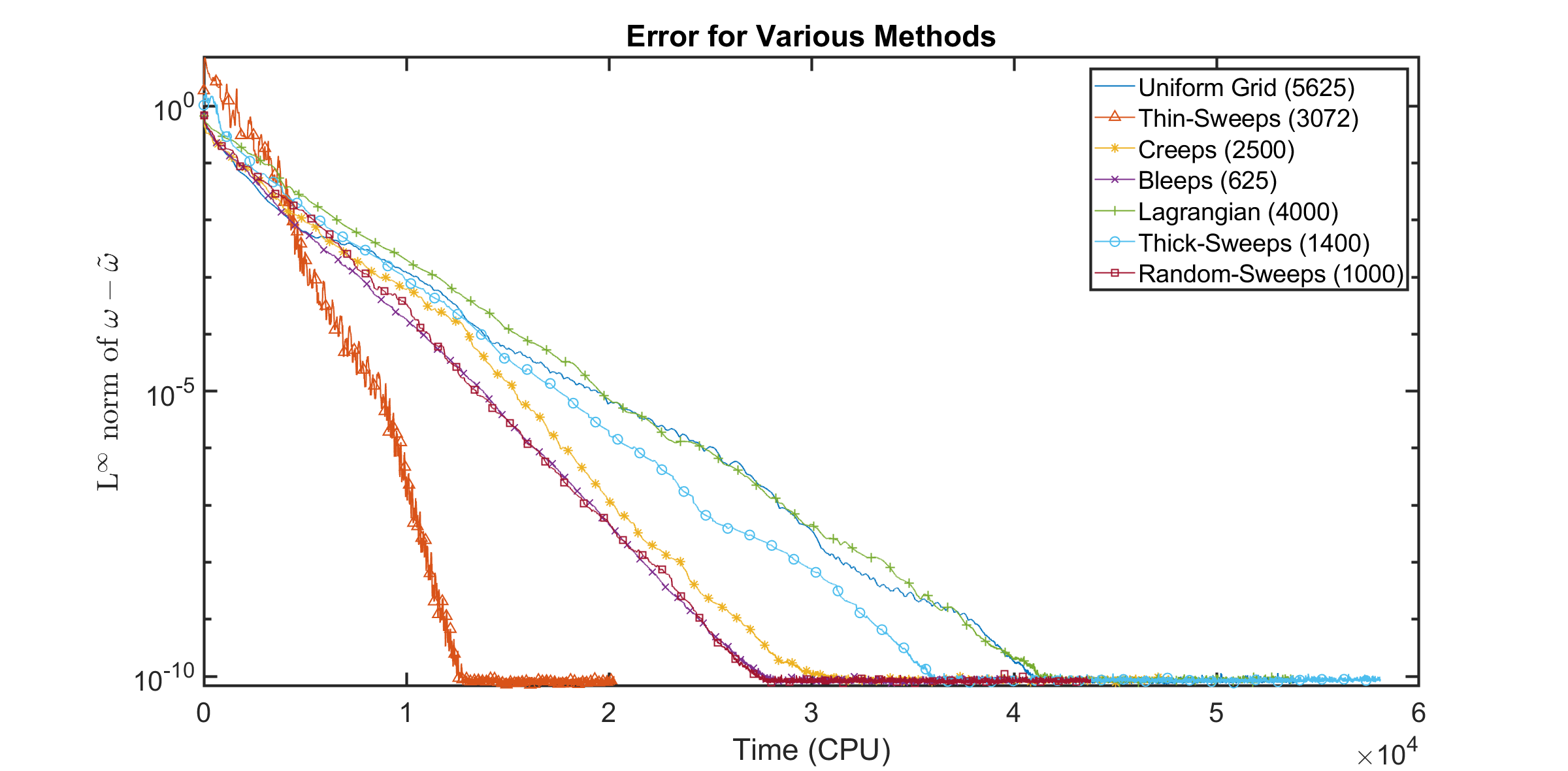}
	\caption{Comparison of error over CPU time for static uniform grid and dynamic observer strategies with number of observers optimized to ensure convergence around the time $t = 70$ (log-linear plot).
	See legend entries for the number of observers used in each method. $\mu = 10$ for all methods except thin-sweeps, which uses $\mu = 30$. Thick-sweeps has constant $x$ velocity $3\frac{\Delta x}{\Delta t}$.
	}
	\label{fig:CPU_min_omega2}
\end{figure}

\begin{figure}
	\includegraphics[width=\textwidth]{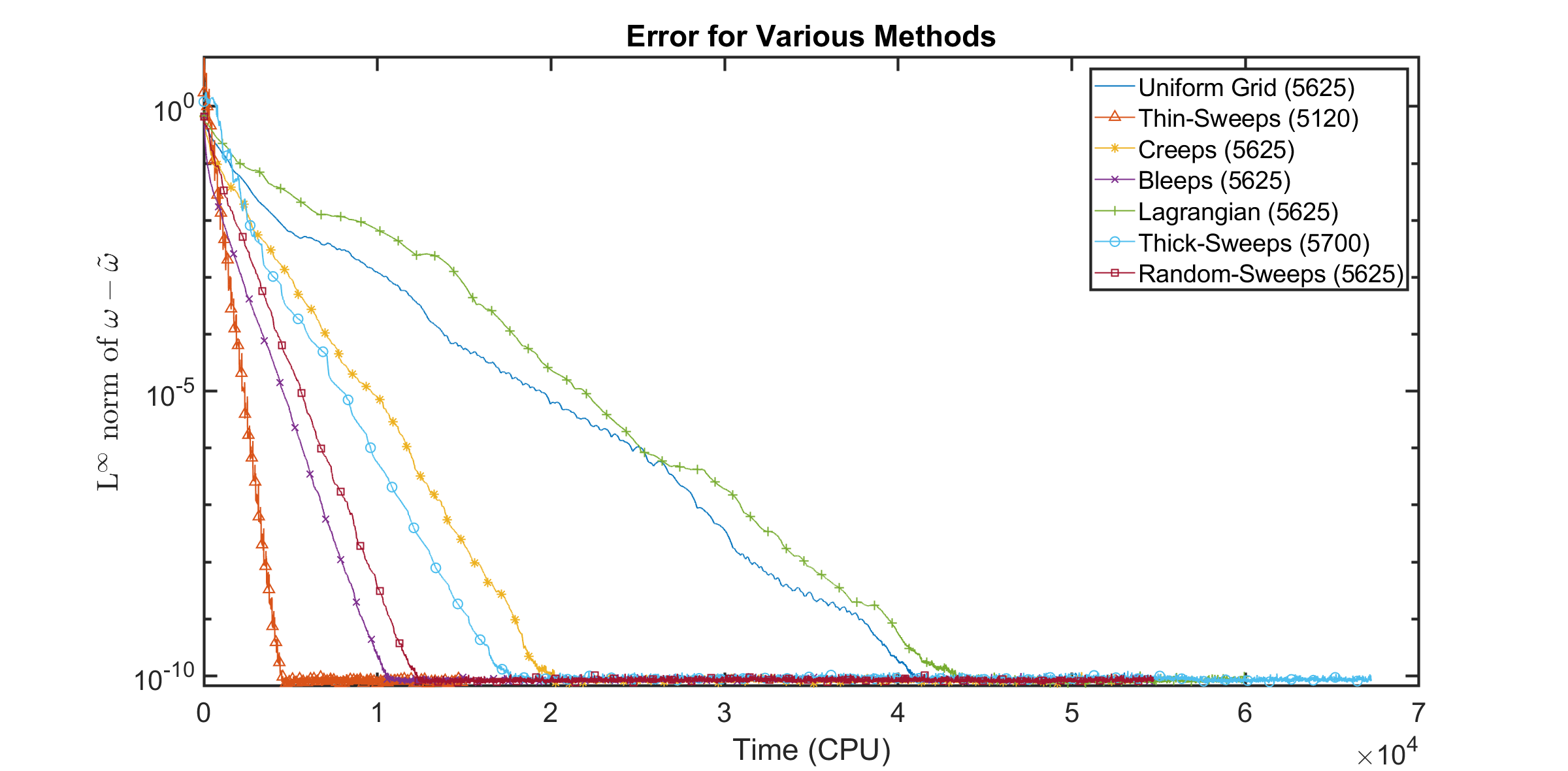}
	\caption{Comparison of error over CPU time for static uniform grid and dynamic observer strategies all initialized using approximately 5,625 observers (log-linear plot). 
	See legend entries for the number of observers used in each method. $\mu = 10$ for all methods except thin-sweeps, which uses $\mu = 30$. Thick-sweeps has constant $x$ velocity $3\frac{\Delta x}{\Delta t}$.
	}	\label{fig:CPU_equal_omega2}
\end{figure}

\begin{figure}
	\includegraphics[width=\textwidth]{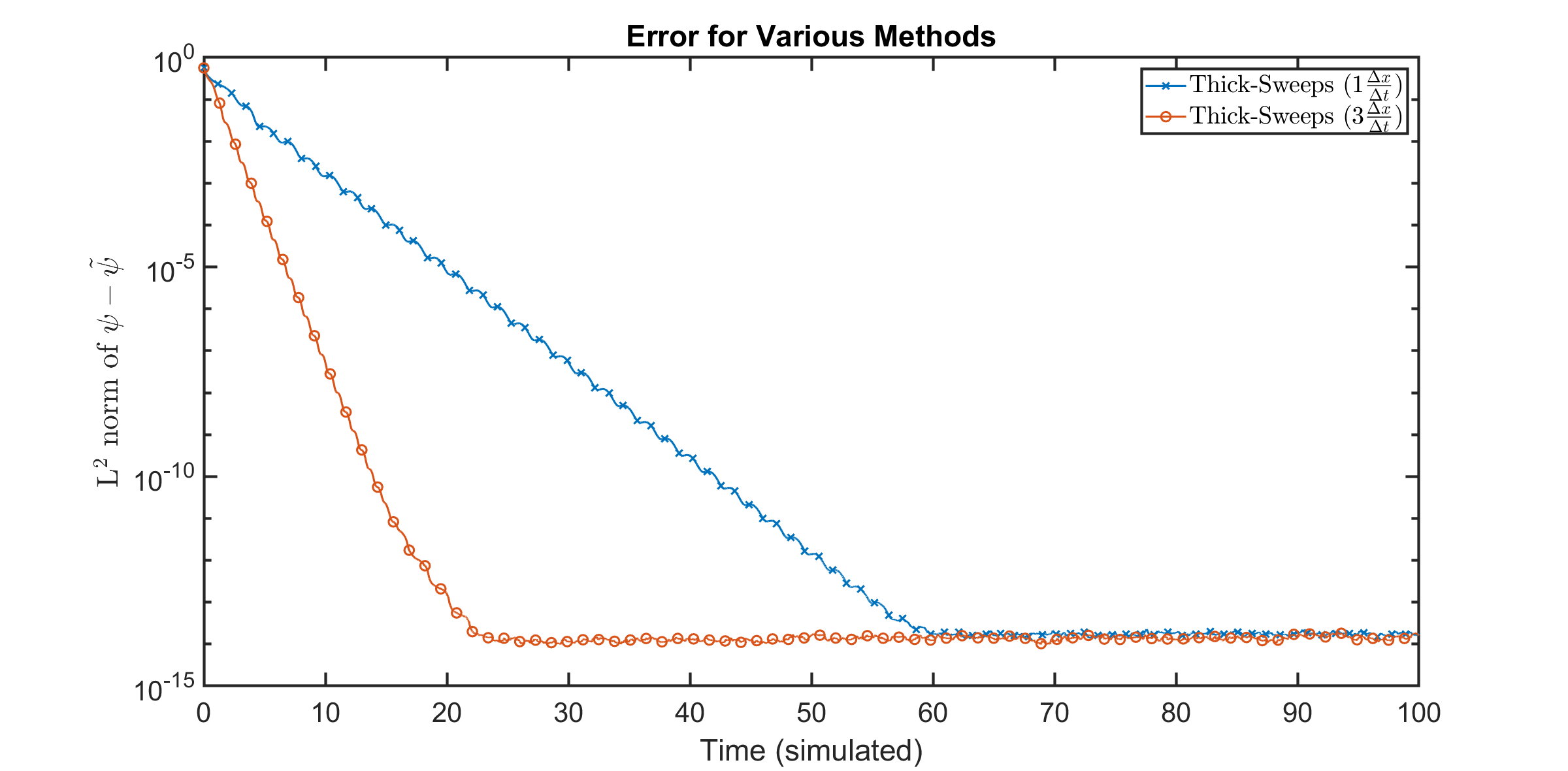}
	\caption{Comparison of Thick-Sweeps methods initialized with different speeds (log-linear plot). Both trials were initialized with the 5,700 observers, $\mu = 10$, moving with constant velocity in the $x$ direction as given in legend parentheses. 
	}
	\label{fig:SIM_thick_comparison}
\end{figure}

\begin{figure}
	\includegraphics[width=\textwidth]{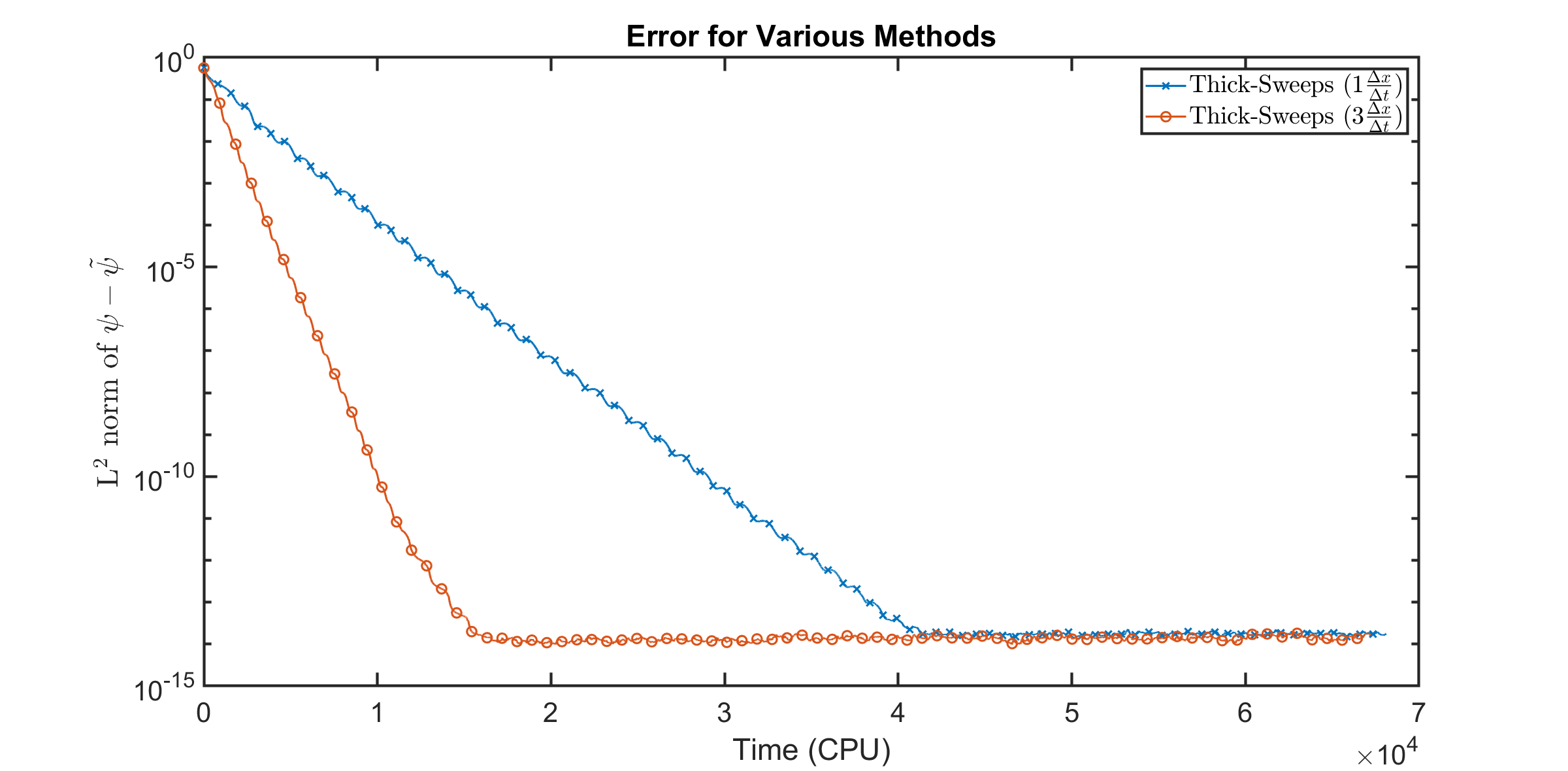}
	\caption{Comparison of Thick-Sweep methods initialized with different speeds (log-linear plot).. Both trials were initialized with the 5,700 observers, $\mu = 10$, moving with constant velocity in the $x$ direction as given in legend parentheses.  
	}
	\label{fig:CPU_thick_comparison}
\end{figure}

\FloatBarrier

\section*{Acknowledgments}
 \noindent 
Author T.E.F. acknowledges the financial support of the USDA National Institute of Food and Agriculture Hatch project \#1020768 and project \#2019-67021-29312. 
Author A.L. acknowledges the financial support of National Science Foundation (NSF) grants DMS-1716801 and CMMI-1953346.
Author C.V. acknowledges the financial support of NSF GRFP grant DMS-1610400.  



\end{document}